%&amstex
\documentstyle{amsppt}
\magnification=1200
\NoRunningHeads
\TagsOnRight

\define\la{\lambda}
\define\La{\Lambda}
\define\R{{\Bbb R}}
\define\C{{\Bbb C}}
\define\Y{{\Bbb Y}}
\define\Z{{\Bbb Z}}
\define\K{{\Bbb K}}
\define\wtp{\widetilde p}
\define\wts{\widetilde s}
\define\wtf{\widetilde f}
\define\wth{\widetilde h}
\define\wtm{\widetilde m}
\define\wtpi{\widetilde\pi}
\define\wtPi{\widetilde\Pi}
\define\Om{\Omega}
\define\om{\omega}
\define\al{\alpha}
\define\be{\beta}
\define\ga{\gamma}
\define\Ga{\Gamma}
\define\de{\delta}
\define\De{\Delta}
\define\ka{\kappa}
\define\si{\sigma}
\define\ep{\varepsilon}
\define\ze{\zeta}
\define\Pz{P_{zz'}}
\define\PP{{\Cal P}}
\define\PPz{\PP_{zz'}}
\define\PDt{\Cal{PD}(t)}
\define\siz{\sigma^{(zz')}}
\define\sit{\sigma^{(t)}}
\define\whsit{\widehat\sigma^{(t)}}
\define\rz{\rho^{(zz')}}
\define\rt{\rho^{(t)}}
\define\tht{\thetag}
\define\M{M_{zz'}}
\define\nuom{\nu^{(\om)}}
\define\Fp{{\Cal F}^+}
\define\Fm{{\Cal F}^-}
\define\Fpm{{\Cal F}^\pm}

\topmatter
\title Point processes and \\
the infinite
symmetric group. \\
Part I: The general formalism and\\
the density function
\endtitle
\author Grigori Olshanski \endauthor
\address{\rm  Dobrushin Mathematics Laboratory, 
Institute for Problems of Information Transmission, 
Bolshoy Karetny 19, 101447 Moscow GSP--4, Russia.} {\it
E-mail\/}: {\tt olsh\@ippi.ras.ru, olsh\@glasnet.ru } 
\endaddress

\thanks Supported by the Russian Foundation for Basic Research 
under grant 98-01-00303 and by
the Russian Program for Support of Scientific Schools under grant
96-15-96060. 
\endthanks

\abstract
We study a 2-parametric family of probability measures on an
infinite--dimensional simplex (the Thoma simplex). These measures
originate in harmonic analysis on the infinite symmetric group
(S.~Kerov, G.~Olshanski and A.~Vershik, Comptes Rendus Acad. Sci.
Paris I 316 (1993), 773-778). Our approach is to interprete them as
probability distributions on a space of point configurations, i.e., as
certain point stochastic processes, and to find the correlation
functions of these processes.  

In the present paper we relate the correlation functions to the
solutions of certain multidimensional moment problems. Then we
calculate the first correlation function which leads to a conclusion
about the support of the initial measures. In the appendix, we discuss a
parallel but more elementary theory related to the well--known
Poisson--Dirichlet distribution.

The higher correlation functions are explicitly calculated in the
subsequent paper (A.~Borodin). In the third part (A.~Borodin and
G.~Olshanski) we discuss some applications and relationships with the
random matrix theory.  

The goal of our work is to understand new phenomena in noncommutative
harmonic analysis which arise when the irreducible representations
depend on countably many continuous parameters.
\endabstract
\toc
\widestnumber\head{\S1.}
\head {} Introduction\endhead
\head \S1. Coherent systems of distributions on the Young graph
\endhead 
\head \S2. The coherent $z$-systems \endhead
\head \S3. Controlling measures \endhead
\head \S4. Point processes \endhead
\head \S5. The density function \endhead
\head \S6. An application \endhead
\head \S7. Appendix: Correlation functions of the Poisson--Dirichlet
\newline processes \endhead
\head \S8. Appendix (A.~Borodin): A proof of Theorem 2.1 \endhead
\head {} References \endhead
\endtoc

\endtopmatter

\document

In this paper, we begin study of a remarkable family of stochastic
point processes. These processes live on the punctured interval
$[-1,1]\setminus\{0\}$ and depend on two real parameters. Our purpose
is to calculate their correlation functions which supply substantial
information about the processes. The present paper is the first one in
a series of papers. It contains introduction to the subject,
description of the method, calculation of the first correlation
function (also called the density function) and an application
concerning the `support' of the processes. The higher correlation
functions are calculated in the subsequent paper \cite{B} by Alexei
Borodin. In the third paper \cite{BO} we discuss certain applications.

The point processes in question originated from harmonic analysis on
the infinite symmetric group \cite{KOV}: they govern decomposition
of the so-called generalized regular representations. I shall briefly
discuss the link with representation theory, as this is the main
motivation of the work. I believe that this new kind of a relationship
between representations and probabilities is interesting. However, in
the body of the paper, we are not dealing with representations, and I
tried to make the exposition formally independent of \cite{KOV}
and accessible to non--experts in representation theory.

Starting with the infinite symmetric group $S(\infty)=\varinjlim S(n)$
(the union of the finite symmetric groups), we form a `$(G,K)$-pair',
where $G$ is the product $S(\infty)\times S(\infty)$ and $K$ is the
diagonal subgroup in $G$.  (Let us emphasize at once that irreducible
representations of $G$, generally speaking, are {\it not\/} tensor
products of two irredicible representations of $S(\infty)$ as it would
be for a `tame' group in place of $S(\infty)$; but $S(\infty)$ is {\it
not tame\/}.) A {\it spherical representation\/} of $(G,K)$ is a
couple $(T,v)$, where $T$ is a unitary representation of $G$ and $v$
is a distinguished $K$-invariant unit vector in the space of $T$. Note
that if $T$ is irreducible then $v$ is unique up to a scalar
factor. Note also that the spherical representations are a particular
case of more general `admissible representations' of $(G,K)$, see
\cite{O1, Ok}.

There exists a parametrization 
$\omega\leftrightarrow (T^{(\omega)},v^{(\omega)})$ 
of irreducible spherical representations by points
$\omega$ of an infinite--dimensional simplex $\Omega$; the latter is
called the {\it Thoma simplex\/}, see \cite{T, VK, O1, Ok}. Any cyclic
representation $(T, v)$ can be decomposed into a direct integral,
$$
T=\int_\Omega T^{(\omega)}P(d\omega),\quad 
v=\int_\Omega v^{(\omega)}P(d\omega). \tag0.1
$$
Here $P$ is a probability measure on $\Omega$, which is uniquely
defined; it is called the {\it spectral measure\/} for $(T,v)$.

In \cite{KOV}, we constructed a family $\{T_z\}$ of admissible
representations of $(G,K)$ depending on a complex parameter $z$. Each
$T_z$ is realized in a $L^2$ space on a compactification $X$ of the
discrete space $S(\infty)$. When $z$ tends to infinity, $T_z$
approaches the conventional biregular representation $T_\infty$ of $G$
in the $\ell^2$ space on $S(\infty)$, so that the representations
$T_z$ form a deformation of $T_\infty$; we call them the {\it
generalized regular representations\/}.

As is well known, the representation $T_\infty$ is irreducible.
However, the representations $T_z$ are highly reducible (as ``true''
regular representations should be). Their decomposition can be viewed
as a model problem of noncommutative harmonic analysis in the
situation when the irreducible representations depend on infinitely
many parameters.  

The construction of \cite{KOV} shows that each $T_z$ possesses a
distinguished $K$-invariant vector $v$; moreover, $v$ is cyclic
provided that $z$ is not integral. Applying the abstract decomposition
\thetag{0.1}, we arrive to a family $\{P_z\}$ of spectral measures on
the simplex $\Omega$. The probability measures $P_z$ are the main
object of the paper. \footnote{Actually, we also consider a
``complementary series'' of spectral measures. In the text, we use the
notation $P_{zz'}$, where either $z'=\bar z$ (the ``principal
series'') or both $z$ and $z'$ are real and satusfy some extra
condition (the ``complementary series'').} Note that they are pairwise
disjoint (\cite{KOV, Theorem 5.3}).

To study the measures $P_z$ we propose the following approach:

1) We define an embedding of the Thoma simplex $\Omega$ to the space
$\Xi$ of configurations in the locally compact space 
$I=[-1,1]\setminus\{0\}$ (a configuration in $I$ is an unordered
collection of points which can accumulate only at 0). Then any
probability measure on $\Omega$ will define a {\it random\/}
configuration in $I$, i.e., a stochastic point process. Thus, the
measures $P_z$ can be interpreted as certain point processes
$\PP_z$.

2) As a characteristic of the point processes $\PP_z$ we
choose the correlation functions. Let $n=1,2,\dots$ and $x_1,\dots,
x_n$ be an arbitrary $n$-tuple of points in $I$. Given a point
process, the probability of the event that a random configuration
intersects each of the infinitely small intervals
$x_1+dx_1,\,\dots,\,x_n+dx_n$ has the form 
$\rho_n(x_1,\dots,x_n)dx_1\dots dx_n$, and the density $\rho_n$ is
called the $n$th {\it correlation function\/}. 
All point processes originated from probability measures on
$\Omega$ are uniquely determined by their correlation functions.

3) We show that the correlation functions of $\PP_z$ can be
obtained from a sequence $\sigma_1,\sigma_2,\dots$ of probability
measures called the {\it controlling measures\/}. The $n$th
controlling measure lives on the $n$-dimensional cube $[-1,1]^n$, and
one can write down all its moments. This reduces the problem of
calculating the functions to a certain multidimensional
moment problem.  

The paper is organized as follows.

\S1 contains preliminaries on symmetric fumctions, the Thoma simplex,
coherent systems of distributions on the Young graph, and their
spectral measures. Using this formalism, we introduce in \S2 the
so--called coherent $z$-systems on the Young graph, by means of which
we define the spectral measures $P_z$. In \S3 we discuss the controlling
measures and their moments, and in \S4 we pass to the point processes.

The technical part of the work begins in \S5. Here we solve a moment
problem and calculate the density function of $\PP_z$. 
We present an integral representation of the density function, Theorem
5.8, and an explicit expression in terms of a multivariate
hypergeometric function (the Lauricella function $F_B$ in three
variables), Theorem 5.12. 

In \S6 we show that the first controlling measure has no atom at 0,
which implies that the measures $P_z$ are concentrated on a
distinguished face of the simplex $\Omega$ (this result was
announced in \cite{KOV}).

There are two appendices.

In the first appendix (\S7) we discuss the Poisson--Dirichlet distributions
$PD(t)$. These distributions were intensively studied in literature 
and they play an important role in the construction \cite{KOV} of the
representations $T_z$. Our purpose is to derive the Watterson \cite{W}
formula for the correlation functions of $PD(t)$ by employing the
general formalism of \S\S 1--4. 

The second appendix (\S8) contains a simple direct proof, due to
A.~Borodin, of Theorem 2.1 asserting the existence of the coherent
$z$-systems.

The results of the present paper were obtained, in the main, in 1992.
Then Alexei Borodin succeeded to calculate
the higher correlation functions; his results constitute the second
part of the work, see \cite{B}. I am very grateful to him for numerous
discussions which exerted a substantial influence on the final version
of the paper. In particular, one of the devices of \cite{B} allowed me
to simplify the derivation of Theorem 5.2.

\head \S1. Coherent systems of distribitions on the Young graph
\endhead

\subhead Symmetric functions \cite{M} \endsubhead
Let $\La$ denote the algebra of symmetric functions over the base
field $\R$. Formally, $\La$ may be defined as $\R[p_1,p_2,\dots]$,
the algebra of polynomials over infinitely many indeterminates
$p_1,p_2,\dots$, called the {\it power sums\/}. Another system of
generators of $\La$ is formed by the {\it complete symmetric
functions\/} $h_1,h_2,\dots$, which are expressed through
$p_1,p_2,\dots$ via the the following relation of generating series
with formal parameter $t$:
$$
1+\sum_{n\ge1} h_nt=\exp(\sum_{n\ge1}(p_n/n)t^n).
$$

A {\it partition\/} is a weakly decreasing sequence
$\la=(\la_1,\la_2,\dots)$ of nonnegative integers with finitely many
nonzero terms. The number of nonzero terms is called the {\it
length\/} of $\la$ and denoted as $\ell(\la)$. Partitions are also
viewed as Young diagrams. By $|\la|$ we denote the sum
$\la_1+\la_2+\dots$ or, equivalently, the number of boxes in the
diagram $\la$. The zero partition (or the empty Young diagram) is
denoted as $\varnothing$.

The elements 
$$
p_\rho=p_{\rho_1}p_{\rho_2}\dots p_{\rho_l},\quad l=\ell(\rho), 
$$
where $\rho$ ranges over the set of partitions, form a basis in $\La$. 
Another distinguished basis in $\La$ is formed by the {\it Schur
functions\/}, which are also indexed by arbitrary partitions and can
be defined, in terms of $h_1,h_2,\dots$, by the {\it Jacobi--Trudi
formula\/}
$$
s_\la=\det[h_{\la_i-i+j}], 
$$
where the order of determinant is $\ell(\la)$ and it is assumed that
$h_0=1$ and $h_n=0$ when $n<0$.

We shall need two important formulas:
$$
\gather
s_\mu\cdot  p_1=\sum_{\la:\;\la\searrow\mu} s_\la,\tag1.1\\
p_\rho=\sum_{\la:\; |\la|=|\rho|} \chi^\la_\rho\, s_\la, \tag1.2
\endgather
$$
see \cite{M, I.5.16 and I.7.8}. Here and in what follows the notation
$\la\searrow\mu$ (or, equivalently, $\mu\nearrow\la$) means that the
diagram $\la$ contains the diagram $\mu$ and differs from it by a
single box (in particular, we have $|\la|=|\mu|+1$). By $\chi^\la$ we
denote the irreducible character indexed by $\la$ (it is a character
of the symmetric group of degree $|\la|$), and $\chi^\la_\rho$ is the
value of $\chi^\la$ on the conjugacy class indexed by $\rho$, see
\cite{M, \S{I}.7}.

In the customary realization of the algebra $\La$, its generators
$p_n$ are identified with the expressions $x_1^n+x_2^n+\ldots$ in
indeterminates $x_1,x_2,\dots$, so that elements of $\La$ become 
symmetric functions in $x_1,x_2,\dots$. But we shall need
another realization, obtained by specializing the generators
$p_n\in\La$ to the following expressions in the indeterminates
$\alpha=(\alpha_1,\alpha_2,\dots)$, $\beta=(\beta_1,\beta_2,\dots)$,
and $\gamma$:
$$
\gathered
p_1\;\mapsto\;\wtp_1(\alpha,\beta,\gamma):=
\sum_{i\ge1}\alpha_i+\sum_{i\ge1}\beta_i+\gamma\\
p_n\;\mapsto\;\wtp_n(\alpha,\beta,\gamma):=
\sum_{i\ge1}\alpha_i^n+(-1)^{n-1}\sum_{i\ge1}\beta_i^n, 
\quad n\ge2, 
\endgathered
\tag1.3
$$
which is equivalent to
$$
1+\sum_{n\ge1}h_nt^n\;\mapsto\; e^{\gamma t}
\prod_{i\ge1}\frac{1+\beta_i t}{1-\alpha_i t}.
$$
This is a generalization of the well--known ``super'' realization of
$\La$; indeed, setting $\gamma=0$ converts the above expressions to
``supersymmetric'' functions in $\alpha$ and $-\beta$, see \cite{S}
and \cite{M, \S{I}.3, Ex. 23}.

\subhead The Thoma simplex \cite{VK, KV, KOO} \endsubhead
We shall abbreviate
$$
\alpha=(\alpha_1,\alpha_2,\dots),\quad
\beta=(\beta_1,\beta_2,\dots).
$$
Let $\Om$ be the set of the triples $\om=(\al,\be,\ga)$ such that
$$
\gathered
\al_1\ge\al_2\ge\dots\ge0,\quad
\be_1\ge\be_2\ge\dots\ge0,\quad
\ga\ge0,\\
\sum_{i\ge1}(\al_i+\be_i)+\ga=1.
\endgathered
\tag1.4
$$
Since $\ga$ is determined by $\al$ and $\be$, we shall sometimes omit
$\ga$ and write $\om=(\al,\be)$. 

The set $\Om$ is an infinite--dimensional simplex;
its faces of codimension 1 have the form
$$
\gather
\Om_i=\{\om\in\Om\bigm|\al_i=\al_{i+1}\},\quad
\Om_{-i}=\{\om\in\Om\bigm|\be_i=\be_{i+1}\},\qquad i\ge1,
\tag1.5\\
\Om_0=\{\om\in\Om\bigm|\ga=0\}. \tag1.6
\endgather
$$

The simplex $\Om$ is called the {\it Thoma simplex\/} in connection
with the pioneering Thoma's work \cite{T}. We equip $\Om$ with the
weakest topology in which the coordinates $\al_i$ and $\be_i$ (but
not $\ga$) are continuous functions. In this topology, $\Om$ is a
metrizable compact space, and the face $\Om_0$ is a dense subset.

We define the functions $\wtp(\om)$ on $\Om$ by setting
$\om=(\al,\beta,\ga)\in\Om$ in the formulas \tht{1.3}. Note that
$\wtp_1(\om)\equiv1$. It is readily verified (\cite{KOO, Lemma 5.2})
that these functions are continuous on $\Om$. \footnote{Observe that
the expression $\sum(\al_i+\be_i)$ is not continuous on $\Om$.}
Consequently, any element $f\in\La$ can be converted to a continuous
function $\wtf$ on $\Om$ by writing $f$ as a polynomial in the
generators $p_n$ and replacing then each $p_n$ by $\wtp_n$. We shall
call $\wtf$ the {\it extended version\/} of $f$, cf.
\cite{KV, KOO}. In particular, we shall deal
with the {\it extended power sums\/} $\wtp_n$ and the {\it extended
Schur functions\/} $\wts_\la$.  

Note that the functions $\wtf$, where $f$ ranges over $\La$, form a
dense subalgebra in the Banach algebra $C(\Om)$ of (real) continuous 
functions on the compact space $\Om$, see \cite{KOO, Lemma 5.3}. 

The algebra $\La$ possesses an involutive automorphism (denoted as
$\om$ in \cite{M}) such that
$$
p_n\,\mapsto\, (-1)^{n-1}p_n,\quad s_\la\,\mapsto\,s_{\la^t}\,,
$$
where $\la^t$ means the transposed diagram. In terms of the
realization \tht{1.3} this involution exactly corresponds to the
symmetry $\al\leftrightarrow\be$. The latter symmetry also defines a
symmetry of the Thoma simplex, which will be denoted as 
$\om\mapsto\om^t$. It follows that
$$
\wts_{\la^t}(\om)=\wts_\la(\om^t).
$$

\subhead Harmonic functions and coherent systems of distributions on
the Young graph \cite{V, VK, KV, K, KOO} \endsubhead
By definition, the vertices of the {\it Young graph\/} $\Y$ are
arbitrary Young diagrams (including $\varnothing$), and its
(oriented) edges are formed by the couples $\mu, \la$ such that
$\mu\nearrow\la$. The number of oriented paths from $\varnothing$ to
$\la$ is called the {\it dimension\/} of $\la$ and denoted as
$\dim\la$ (we agree that $\dim\varnothing=1$). The function $\dim$ on
the vertices of $\Y$ satisfies the {\it recurrence relation} 
$$
\dim\la=\sum_{\mu:\; \mu\nearrow\la} \dim\mu.
$$

These definitions are inspired by the {\it Young branching rule\/}
for the characters of the finite symmetric groups $S(n)$,
$n=1,2,\dots$, 
$$
\chi^\la\bigm|_{S(n-1)}=\sum_{\mu:\; \mu\nearrow\la} \chi^\mu,
\quad n=|\la|.
$$
(See \cite{V, JK, OV}.) The above two relations
show that the dimension of a diagram $\la$ coincides with the
dimension of the character $\chi^\la$, i.e., with the number
$\chi^\la_{(1^n)}$, $n=|\la|$.  

A (real) function $\varphi(\la)$ on the vertices of $\Y$ is called
{\it harmonic\/} if it satisfies the ``harmonicity condition''
$$
\varphi(\mu)=\sum_{\la:\;\la\searrow\mu} \varphi(\la)
$$
for each diagram $\mu$. Let $\Y_n$ denote the set of Young diagrams
with $n$ boxes (equivalently, the set of partitions of the number
$n$). The sets $\Y_n$ define a {\it grading\/} of the graph $\Y$. 
Clearly, the knowledge of a harmonic function on a ``floor'' $\Y_n$
determines it on the preceding floors $\Y_{n-1},\Y_{n-2},\dots$.

Let $M(\la)$ be a function on the vertices of $\Y$ and $M_n$ denote
its restriction to the $n$th floor $\Y_n$, $n=0,1,2,\dots$. We call
$M=(M_n)$ a {\it coherent system of distributions\/} on the Young
graph (coherent system, for short) if $M$ is nonnegative, normalized
at $\la=\varnothing$ (i.e., $M(\varnothing)=M_0(\varnothing)=1$), and
the function $M(\la)/\dim\la$ is harmonic, i.e.,
$$
M(\mu)=\sum_{\la:\;\la\searrow\mu} \frac{\dim\mu}{\dim\la}\, M(\la).
$$
The latter formula and the recurrence relation for the dimension
function imply that the $M_n$'s are probability distributions on the
floors of the Young graph (see \cite{KOO, Lemma 8.1}). 

By a measure on $\Om$ we shall always mean a Borel measure with
respect to the canonical Borel structure of the topological space
$\Om$. According to \cite{KOO, Theorem B}, there is a bijective
correspondence $M\leftrightarrow P$ between the coherent systems $M$ on
the Young graph and the probability measures $P$ on the Thoma
simplex. This correspondence is characterized by the relation
$$
\frac{M(\la)}{\dim\la}=\int_\Om \wts_\la(\om)\,P(d\om), \tag1.7
$$
where $\la$ is an arbitrary Young diagram and $\wts_\la$ is the
extended Schur function as defined above. Moreover, as is shown in
the proof of this result, the measure $P$ is approached, in a certain
sense, by the finite probability distributions $M_n$ as
$n\to\infty$. We shall call $P$ the {\it spectral measure\/} of $M$.

If $M$ is a coherent system on $\Y$ then the function 
$M^t(\la)=M(\la^t)$ is a coherent system, too (indeed, this follows from
the fact that the involution $\la\mapsto\la^t$ is a symmetry of the
Young graph). On the other hand, if $P$ is a measure on $\Om$, let
$P^t$ denote the image of $P$ under the symmetry $\om\mapsto\om^t$
defined above. Now, we have
$$
\text{if}\quad M\,\leftrightarrow\,P\quad 
\text{then}\quad M^t\,\leftrightarrow\,P^t. \tag1.8
$$ 
Indeed, this follows from \tht{1.7} and the equality
$\wts_{\la^t}(\om)=\wts_\la(\om^t)$ mentioned above. 

\subhead Connection with representations \cite{T, VK, KV, O1, O2}
\endsubhead 
The constructions described above are inspired by the representation
theory of the infinite symmetric group $S(\infty)=\varinjlim S(n)$.
Let $M$ be a coherent system on $\Y$ and $\varphi(\la)=M(\la)/\dim\la$
the corresponding harmonic function. For each $n$, the linear
combination of irreducible characters
$$
\chi_n=\sum_{\la:\;|\la|=n}\varphi(\la)\chi^\la
$$
is a central positive definite function on $S(n)$, normalized at the
unit element. By the harmonicity condition, the functions
$\chi_1,\chi_2,\dots$ are compatible with the embeddings
$S(n-1)\hookrightarrow S(n)$. Consequently, they define a central, 
positive definite, normalized function $\chi$ on the group
$S(\infty)$. Let $X$ denote the set of all such functions. This is a
convex set; its extreme points are called the {\it characters\/} of
the group $S(\infty)$ (in the sense of von Neumann). 

Via the Gelfand--Naimark--Segal construction, characters generate
{\it finite factor representations\/} of the group $S(\infty)$. They
also correspond to {\it irreducible unitary spherical
representations\/} of the Gelfand pair $(G,K)$, where $G$ 
stands for the ``bisymmetric group'' $S(\infty)\times S(\infty)$ and
$K$ is the diagonal subgroup of $G$, see \cite{O1, O2}. 

The integral representation \tht{1.7} has the following meaning. First
of all, via the correspondence $M\leftrightarrow\chi$, \tht{1.7}
implies that the characters are parametrized by the points
$\om\in\Om$ (we shall write them as $\chi^{(\om)}$). \footnote{Thus,
the Thoma simplex may be viewed as the {\it spherical dual\/} of the
Gelfand pair $(G,K)$.} This result is known as {\it Thoma's
theorem\/} \cite{T}; see also \cite{VK}. The extreme
coherent system $M^{(\om)}\leftrightarrow\chi^{(\om)}$ is given by the
formula 
$$
M^{(\om)}(\la)=\dim\la\cdot \wts_\la(\om),
$$
which is equivalent to {\it Thoma's formula} \cite{T} 
$$
\chi^{(\om)}(\rho)=\wtp_{\rho_1}(\om)\wtp_{\rho_2}(\om)\dots\ .
$$
In the latter expression, the left--hand side is the value of
$\chi^{(\om)}$ at the conjugacy class 
in $S(\infty)$ of an arbitrary cycle--type 
$\rho=(\rho_1,\rho_2,\dots,1,1,\dots)$, and the right--hand side is
correctly defined, because almost all cycles have length 1 and
$\wtp_1\equiv1$. 

Next, \tht{1.7} implies that any function $\chi\in X$ is uniquely
decomposed into a continual convex combination of the characters
$\chi^{(\om)}$. As $\chi$ generates a cyclic spherical representation
$T$ of the pair $(G,K)$, this also means that the spectral measure
$P$ effectues the decomposition of $T$ in a direct integral of
irreducible spherical representations. 

\head \S2. The coherent $z$-systems \endhead

Set
$$
\M(\la)=\frac{\prod\limits_{(ij)\in\la}(z+j-i)(z'+j-i)}
{(t)_{|\la|}}\cdot 
\frac{\dim^2\la}{|\la|\ !}\,.\tag2.1
$$
Here $\la$ is an arbitrary Young diagram, $z$ and $z'$ are complex
parameters, $(ij)\in\la$ is an arbitrary box of $\la$ ($i$ and $j$
are the numbers of the row and the column containing the box),
$t=zz'$, and an expression of type $(t)_n$ means the Pohgammer
symbol: 
$$
(t)_n=\Gamma(t+n)/\Gamma(t)=t(t+1)\dots(t+n-1);
$$
\tht{2.1} is correctly defined if $t\ne0,-1,-2,\dots$. We also agree
that $\M(\varnothing)=1$.  

\proclaim{Theorem 2.1} The function $\M(\la)/\dim\la$ is harmonic. 
\endproclaim

\demo{Proof} This follows from \cite{KOV, Theorem 3.1}. A direct
combinatorial argument is given in Kerov's paper \cite{Ke2} (actually,
\cite{Ke2} contains a more general result). Another combinatorial
proof was proposed by Postnikov \cite{P}. In \S8 below we present one more
direct proof, due to Borodin. \qed
\enddemo

\proclaim{Proposition 2.2} Assume that $z,z'$ are not integers and
are such that $t=zz'$ is not equal to $0,-1,-2,\dots$. Then
$\M(\la)\ne0$ for any $\la$. 
Moreover, $\M$ is strictly positive if and only if the parameters
$z,z'$ satisfy one of the following conditions:

{\rm(i)} $z'=\bar z$ and $z\notin\Z${\rm;}

{\rm(ii)} $z$ and $z'$ are real and are both contained in an 
open interval of the form $(m,m+1)$ with $m\in\Z$.
\endproclaim

\demo{Proof} The first claim is obvious from \tht{2.1}, let us
check the second claim. Assume that (i) or (ii) holds. Then $t>0$ and
$(z+k)(z'+k)>0$ for any $k\in\Z$, which implies $\M(\la)>0$. 

Conversely, assume that $\M(\la)>0$ for any $\la$. Let
$\mu\nearrow\la$ be an edge of $\Y$,  $n=|\mu|$, $(ij)$ be the
box $\la\setminus\mu$, and $k=j-i$. Comparing $\M(\mu)$ and
$\M(\la)$, we see that 
$$
\frac{(z+k)(z'+k)}{(t+n)}>0.
$$
Clearly, this inequality holds for any numbers $k\in\Z$ and
$n=1,2,\dots$ that correspond to edges of $\Y$, that is to say, for
any $k,n$ such that either $k\ne0$ and $n\ge|k|$, or $k=0$ and
$n\ge3$.  

Now, we fix $k$ and let $n\to\infty$. {}From the above inequality we
conclude that the numerator $(z+k)(z'+k)$ must be real and strictly
positive. Since this holds for any $k\in\Z$, both $zz'$ and $z+z'$
are real. It follows that either $z,z'$ are complex conjugate or
both real. It remains to examine the latter possibility. Using the
fact that $(z+k)(z'+k)$ is not only real but strictly positive, we
get condition (ii) of the proposition. \qed
\enddemo

\proclaim{Corollary 2.3} Under assumptions {\rm(i)} or {\rm(ii)} of
Lemma\/ {\rm 2.2}, $\M$ is a coherent system of distributions on the
Young graph. \qed 
\endproclaim

{}From now on we shall assume that the parameters $z,z'$ satisfy one of
the conditions (i), (ii) of Proposition 2.2 (note that $t\ne0,-1,-2,\dots$
holds automatically then). We shall call the coherent systems $\M$ the
(coherent) {\it $z$-systems\/}. To distinguish between (i) and (ii), we shall
speak about the {\it principal series\/} and {\it complementary
series\/} of $z$-systems, respectively. A motivation for such a
terminology is that the difference $z-z'$ ranges, respectively, over
the imaginary axis and the open interval $(-1,1)$, like the parameters
of the principal or the complementary series for $SL(2,\R)$.

The principal series of $z$-systems first appeared in
\cite{KOV}: in that form we described the spherical functions of
the generalized regular representations of the pair $(G,K)$. The
existence of the complementary series was observed in 1995 by
Borodin. The $z$-systems can be characterized as the only
coherent systems of distributions on the Young graph satisfying a
``multiplicativity condition'', see \cite{R}. 

In the present paper we do not deal with a ``degenerate series'' of
coherent systems which arises when one of the parameters
$z,z'$ is integral. The ''degenerate'' coherent systems live, in
essence, on truncated versions of the Young graph. About them, see
\cite{KOV} and \cite{K}. 

\example{Remark 2.4} The expression \tht{2.1} for $\M$ has two evident
symmetries, each of which has a representation--theoretic meaning.
First, \tht{2.1} does not change under the transposition
$z\leftrightarrow z'$, which leads to certain intertwining operators
for the generalized regular representations. Second, we have
$$
\M(\la^t)=M_{-z,-z'}(\la), \tag2.2
$$
where $\la\mapsto\la^t$ denote transposition of Young diagrams. This
symmetry reflects the well--known fact that, for irreducible
characters $\chi^\la$, transposing $\la$ is equivalent to tensoring
$\chi^\la$ with the one--dimensional sign character.
\endexample

In the next two propositions, we rewrite the expression \tht{2.1} by
making use of two explicit formulas for $\dim\la$.

\proclaim{Proposition 2.5} The expression \tht{2.1} can be written as
$$
\gathered
\M(\la)=\frac{|\la|!}{(t)_{|\la|}}
\prod_{i=1}^l(z-i+1)_{\la_i}(z'-i+1)_{\la_i}\\
\times\frac
{\prod\limits_{1\le i,j\le l}(\la_i-\la_j+j-i)^2}
{\prod\limits_{1\le i\le l}((\la_i+l-i)!)^2}\,,
\endgathered
\tag2.3
$$
where $l\ge\ell(\la)$ may be taken arbitrarily.
\endproclaim

\demo{Proof} The first product is exactly the product over the boxes
$(ij)\in\la$ in formula \tht{2.1}, and the remaining terms come from
the formula
$$
\frac{\dim\la}{|\la|!}=\frac
{\prod\limits_{1\le i,j\le l}(\la_i-\la_j+j-i)}
{\prod\limits_{1\le i\le l}(\la_i+l-i)!}\,, 
\quad l\ge\ell(\la). \tag2.4
$$
It is readily seen that the right--hand side is stable on
$l\ge\ell(\la)$. Hence, it suffices to 
check it for a particular value of $l$. But for $l=|\la|$, this
coincides with the formula of \cite{M, \S{I}.7, Ex. 6}. \qed
\enddemo

The explicit expression \tht{2.3} is not quite satisfactory, because
it does not reflect the symmetry \tht{2.2}. A symmetric expression
can be obtained with the help of the Frobenius notation for Young
diagrams \cite{M, \S{I}.1}:
$$
\la=(p_1,\dots,p_d\,|\, q_1,\dots,q_d), \tag2.5
$$
where $d$ is the length of the diagonal in $\la$, and
$$
p_i=\la_i-i,\quad q_i=(\la^t)_i-i, \qquad 1\le i\le d,
$$
are the {\it Frobenius coordinates\/} of $\la$. Note that
$$
|\la|=|p|+|q|+d,
$$
where
$$
|p|=p_1+\dots+p_d,\quad 
|q|=q_1+\dots+q_d.
$$

\proclaim{Proposition 2.6} In the Frobenius notation \tht{2.5}, the
expression \tht{2.1} can be written as follows
$$
\gathered
\M(\la)=\frac{(|p|+|q|+d)!\,t^d}{(t)_{|p|+|q|+d}}\,
\prod_{i=1}^d\frac
{(z+1)_{p_i}(z'+1)_{p_i}(-z+1)_{q_i}(-z'+1)_{q_i}}
{(p_i!)^2\,(q_i!)^2}\\
\times\prod_{i=1}^d\prod_{j=1}^d(p_i+q_j+1)^{-2}
\prod_{1\le i<j\le d}(p_i-p_j)^2(q_i-q_j)^2\,.
\endgathered
\tag2.6
$$
\endproclaim

\demo{Proof} Given a box $(ij)\in\la$, its {\it hook\/} is defined as
the shape formed by the boxes $(kl)\in\la$ such that either 
$k=i,l\ge j$ or $k>i, l=j$. The total number of boxes in the hook is
called the {\it hook--length\/} and denoted as $h(i,j)$.

Let us represent the shape $\la$ as the union of the diagonal
hooks. Then the contribution of the $k$th diagonal hook
($k=1,\dots,d$) to the product 
$$
\prod_{(ij)\in\la}(z+j-i)(z'+j-i)
$$
is equal to
$$
\gather
(z-q_k)\dots(z-1)z(z+1)\dots(z+p_k)\\
\times(z'-q_k)\dots(z'-1)z'(z'+1)\dots(z'+p_k)\\
=zz'(z+1)_{p_k}(z'+1)_{p_k}(-z+1)_{q_k}(-z'+1)_{q_k}\\
=t(z+1)_{p_k}(z'+1)_{p_k}(-z+1)_{q_k}(-z'+1)_{q_k},
\endgather
$$
which explains the term $t^d$ and the first product in 
\tht{2.6}. The remaining terms in \tht{2.6} come from the following
formula expressing $\dim\la$ in the Frobenius notation:
$$
\frac{\dim\la}{|\la|!}=
\frac
{\prod\limits_{1\le i,j\le d}(p_i-p_j)(q_i-q_j)}
{\prod\limits_{1\le i\le d}\,
\prod\limits_{1\le j\le d}
(p_i+q_j+1)\cdot
\prod\limits_{1\le i\le d}(p_i!\ q_i!)}\,.\tag2.7
$$

To check \tht{2.7} we start with the well--known {\it hook formula\/}
$$
\frac{\dim\la}{|\la|!}=
\prod_{(ij)\in\la}h(i,j)^{-1}, \tag2.8
$$
which, by virtue of formula \tht{2.4}, is equivalent to the
identity 
$$
\prod_{(ij)\in\la}h(i,j)^{-1}=
\frac
{\prod\limits_{1\le i,j\le l}(\la_i-\la_j+j-i)}
{\prod\limits_{1\le i\le l}(\la_i+l-i)!}\,, 
\quad l\ge\ell(\la), \tag2.9
$$
see {M, \S{I}.1, Ex. 1}. 

Now, let us divide the shape $\la$ into three pieces: the square
shape of size $d\times d$, the diagram $\la^+$ formed by the boxes
$(ij)$ with $j>d$, and the diagram $\la^-$ formed by the boxes $(ij)$
with $i>d$. The hook--length of a box $(ij)$ from the square shape is
equal to $p_i+q_j+1$. Consequently, the product of the hook--lengths
over the boxes entering the square shape is equal to the double
product 
$$
\prod\limits_{1\le i\le d}\,
\prod\limits_{1\le j\le d}
(p_i+q_j+1)
$$
in the denominator of \tht{2.7}. To explain the remaining terms in
\tht{2.7}, we express the products of the hook--lengths in
the diagrams $\la^+$ and $\la^-$ via the identity \tht{2.9}, where 
we substitute $l=d$ and $\la=\la^+$ or $\la=(\la^-)^t$. It should also
be noted that for any box in $\la^\pm$, its hook with respect to
$\la$ is the same as the hook with respect to $\la^\pm$. \qed
\enddemo

Formula \tht{2.6} will be used for calculations in \S5 below and in
\cite{B}. 

\head \S3. Controlling measures \endhead

Given a point $\om=(\al,\be,\ga)\in\Om$, we define the corresponding
{\it Thoma measure\/} on $[-1,1]$ as
$$
\nuom=\sum_{i\ge1}\al_i\de(\al_i)+
\sum_{i\ge1}\be_i\de(-\be_i)+\ga\de(0),
$$
where $\de(x)$ stands for the Dirac mass at $x\in[-1,1]$. Clearly,
$\nuom$ is a probability measure.

\proclaim{Proposition 3.1} The moments of $\nuom$ are given by the
formula
$$
\int_{-1}^1x^l\ \nuom(dx)=\wtp_{l+1}(\om),
\qquad l=0,1,2,\dots,
$$
where $\wtp_{l+1}$ are the extended power sums as defined in \S1.
\endproclaim

\demo{Proof} This is a direct consequence of the definition of the
Thoma measure and that of $\wtp_n$. \qed
\enddemo

Let $\operatorname{Prob}[-1,1]$ denote the set of probability Borel
measure on $[-1,1]$; this set has a natural Borel structure
\cite{DVJ}. The map
$$
\Om\to\operatorname{Prob}[-1,1], \qquad
\om\mapsto\nuom,
$$
is Borel--measurable (this is a routine exercise). Therefore,
equiping $\Om$ with a probability measure $P$, we obtain a {\it
random measure\/} on $[-1,1]$. Since the above map in injective, we
may interprete probability measures $P$ on $\Om$ as random measures on
$[-1,1]$.

The next construction looks rather natural from the point of view of
the theory of random measures (see, e.g., \cite{DVJ}) or the
exchangeability theory \cite{A}. We take the infinite product
$$
\nuom_\infty=\nuom\times\nuom\times\dots, \qquad \om\in\Om,
$$
which is a probability measure on the infinite--dimensional cube
$$
[-1,1]^\infty=[-1,1]\times[-1,1]\times\dots, 
$$
and we average $\nuom_\infty$ with respect to a given probability
measure $P$ on $\Om$:
$$
\si=\int_\Om\nuom_\infty\ P(d\om).
$$
In other words, viewing $(\Om,P)$ as a probability space, we take the
expectation of the random measure $\nuom_\infty$. We shall $\si$ the
{\it controlling measure\/} (of infinite order) for the measure $P$.
Clearly, $\si$ is a symmetric probability measure on the
infinite--dimensional cube. Its projection on the 
$n$-dimensional cube $[-1,1]^n$ will be called the $n$th 
{\it controlling measure\/} of $P$ and denoted as $\si_n$:
$$
\si_n=\int_\om\,\underbrace{\nuom\times\dots\times\nuom}_
{\text{$n$ times}}\, P(d\om).
$$

For a Young diagram $\la$, let $d(\la$ be the number of diagonal
boxes. By $\chi^\la_{(r_1,\dots,r_n)}$, where $r_1,\dots,r_n$ are
(non necessarily decreasing) numbers $\ge1$ such that
$r_1+\dots+r_n=|\la|$, we denote the value of the irreducible
character at any permutation with $n$ cycles of length
$r_1,\dots,r_n$. 

\proclaim{Proposition 3.2} Let $M$ be a coherent system on the Young
graph, $P$ be its spectral measure on $\Om$, and $\si_n$ be the $n$th
controlling measure of $P$.

{\rm (i)} The moments of $\si_n$ satisfy the relations
$$
\gather
\int\limits_{[-1,1]^n}x_1^{l_1}\dots x_n^{l_n}\,
\si_n(dx_1\dots dx_n)=
\int_\Om\wtp_{l_1+1}(\om)\dots\wtp_{l_n+1}P(d\om) \tag3.1a\\
=\sum\Sb \la:\; d(\la)\le n\\ |\la|=l_1+\dots+l_n+n\endSb
\chi^\la_{(l_1+1,\dots,l_n+1)}\,
\frac{M(\la)}{\dim\la}\,, \tag3.1b
\endgather
$$
where $l_1,\dots,l_n=0,1,2,\dots$.
\endproclaim

\demo{Proof} By definition of $\si_n$, the left--hand side is the
expectation (with respect to $(\Om,P)$) of the integral
$$
\int\limits_{[-1,1]^n}x_1^{l_1}\dots x_n^{l_n}\,
\nuom(dx_1)\dots\nuom(dx_n).
$$
By Proposition 3.1, this integral is equal to
$$
\wtp_{l_1+1}(\om)\dots\wtp_{l_n+1}(\om).
$$
Integrating over $\om\in\Om$ with respect to $P$, we get \tht{3.1a}.

According to the identity \tht{1.2}, the right--hand side of
\tht{3.1a} is equal to
$$
\int_\Om\sum_{\la:\; |\la|=l_1+\dots+l_n+n}
\chi^\la_{(l_1+1,\dots,l_n+1)}\,\wts_\la(\om)P(d\om).
$$

The Murnaghan--Nakayama rule (\cite{M, \S{I}.7, Ex. 5}) implies that
$\chi^\la_{(l_1+1,\dots,l_n+1)}$ vanishes when $d(\la)>n$, so we may
introduce the supplementary requirement $d(\la)\le n$ into the above
sum.

Finally, transposing integration and summation and applying \tht{1.7}
we get \tht{3.1b}.  

\qed
\enddemo

Thus, when $M$ is known, we can, in principle, find the controlling
measures $\si_n$ from the moment problem \tht{3.1}. The latter has a
unique solution, because the support of $\si_n$ is bounded.

The restriction $d(\la)\le n$ that appears in formula \tht{3.1} will
play a crucial role in what follows. It means that the $n$th
controlling measure is completely determined by the restriction of
$M$ to the set of diagrams $\la$ contained in the $\Gamma$-like shape
formed by the boxes $(ij)$ with $\min(i,j)\le n$. In
particular, for calculating $\si_n$, it suffices to know the values
of $M$ on the hook diagrams. 

Let $\Pz$ be the spectral measure of $\M$; we shall refer to the
family $\{\Pz\}$ as to that of {\it spectral $z$-measures\/}. The
$n$th controlling measure of $\Pz$ will be denoted by $\siz_n$. 

\proclaim{Proposition 3.3} The moments of the measure $\siz_n$ are
given by the formula
$$
\gathered
\int\limits_{[-1,1]^n}x_1^{l_1}\dots x_n^{l_n}\,
\siz_n(dx_1\dots dx_n)\\
=\sum_{d=1}^n\; \sum\Sb 
p_1>\dots>p_d\ge0\\
q_1>\dots>q_d\ge0\\
|p|+|q|+d=|l|+n\endSb
\chi^{(p_1,\dots,p_d\,|\,q_1,\dots,q_d)}_{(l_1+1,\dots,l_n+1)}\\
\times \frac{t^d}{(t)_{|p|+|q|+d}}\;
\prod_{i=1}^d
\frac{(z+1)_{p_i}(z'+1)_{p_i}(-z+1)_{q_i}(-z'+1)_{q_i}}
{p_i!\,q_i!}\\
\times \prod_{i,j=1}^d(p_i+q_j+1)^{-1}\;
\prod_{1\le i,j\le d}(p_i-p_j)(q_i-q_j)\,
\endgathered
\tag3.2
$$
where $l_1,\dots,l_n=0,1,2,\dots$ and
$$
|p|=p_1+\dots+p_d,\quad
|q|=q_1+\dots+q_d,\quad
|l|=l_1+\dots+l_n+n.
$$
\endproclaim

\demo{Proof} This follows from Proposition 3.2 and Proposition 2.6.
\enddemo

We close the section by indicating first applications of the
controlling measures: we shall see that $\si_1$ and $\si_2$ control
the location of the spectral measure $P$ with respect to the faces
\tht{1.5}--\tht{1.6} of the Thoma simplex.

\proclaim{Proposition 3.4} A spectral measure $P$ is concentrated on
the face $\Om_0$ {\rm(}see \tht{1.6}{\rm)} if and only if the first
controlling measure $\si_1$ has no atom at $0$.
\endproclaim

\demo{Proof} Let us regard the parameter $\ga$ as a function
$\ga(\om)$ on $\Om$. This function is lower semicontinuous, hence a
Borel function. By definition of $\si_1$ and of $\nuom$, 
$$
\si(\{0\})=\int_\Om\ga(\om)P(d\om).
$$
Since $\ga(\om)$ is nonnegative, this expression is equal to 0
if and only if $\ga(\om)=0$ almost surely with respect to $(\Om,P)$,
and the latter happens if and only $P$ is concentrated on the face
$\Om_0$. \qed
\enddemo

\proclaim{Proposition 3.5} Let $P$ be a spectral measure and $\si_1$,
$\si_2$ be its first two controlling measures. Consider the subset
$$
\Delta=\{(x,x)\bigm| x\ne0\}\subset [-1,1]^2
$$
and identify it with the punctured interval $I=[-1,1]\setminus\{0\}$.

We always have
$$
\si_2\big|_\Delta(dx)\ \ge\  |x|\cdot\si_1\big|_I(dx),
$$ 
{\rm(}the former measure majorates the latter{\rm)}, and both
measures coincide if and only if, for $P$-almost all points
$\om\in\Om$, there is no repetitions of the type
$\al_i=\al_{i+1}\ne0$ or $\be_i=\be_{i+1}\ne0$.  
\endproclaim

Note that if such repetitions occur with a nonzero probability, then
at least one of the faces $\Om_{\pm i}$, $i\ge1$, is not a
$P$-negligible set. 

\demo{Proof} For any $\om\in\Om$ we have
$$
\gather
(\nuom\times\nuom)\big|_\Delta
=\sum\Sb i,j\\ \al_i=\al_j>0\endSb 
\al_i\al_j\ \de(\al_i)\de(\al_j)
+\sum\Sb i,j\\ \be_i=\be_j>0\endSb 
\be_i\be_j\ \de(-\be_i)\de(-\be_j)\\ 
\ge\,
\sum_{i:\; \al_i>0}\, \al_i^2\ \de(\al_i)+
\sum_{i:\; \be_i>0}\, \be_i^2\ \de(-\be_i)
=|\cdot|\ \nuom\big|_I\,,
\endgather
$$
and equality holds if and only there are no couples $i\ne j$ such that
$\al_i=\al_j>0$ or $\be_i=\be_j>0$.

Integrating over $\om$ with respect to $P$ gives the desired claim.
\qed 
\enddemo

\head \S4. Point processes \endhead

We start with generalities about point processes and their correlation
functions; our main reference here is the book \cite{DVJ}.

Let $X$ be a standard Borel space equipped with a ``bornology''. The
latter means that we know what subsets of $X$ are ``bounded''. We
assume that the family of the bounded subsets is closed under taking
finite union and passage to a subset, and the whole space $X$ can be
represented as countable union of bounded subsets. Bounded Borel
subsets will be called {\it test subsets\/}. 

A {\it configuration\/} in $X$ is a finite or
countable system of points $\xi=(x_1,x_2,\dots)$ in $I$ such that its
intersection with any test set is finite. The word ``system'' is
employed to emphasize that $\xi$ is neither a subset of $X$ (since
repetitions are permitted), nor a sequence (since the points are not
ordered); strictly speaking, $\xi$ is a {\it multiset\/}. However, it
is occasionally convenient to regard $\xi$ as the image of a sequence of
points.

The space of configurations in $X$ will be denoted as $\Xi$. 
For a configuration $\xi\in\Xi$ and a test set $A$, let $|\xi\cap A|$
denote the number of points of $\xi$ (counted with their
multiplicities) that occur in $A$. By the assumption, this number is
always finite. We equip $\Xi$ with the Borel structure generated by
the functions of the form $\xi\mapsto |\xi\cap A|$. 

A {\it point process\/} in $X$ is a (Borel) measurable map from a
probability space (the {\it state space\/}) to $\Xi$. Equivalently, a
point process defines a {\it random configuration\/} in $X$. A point 
process is called {\it simple\/} if the random configuration has no
multiple points almost surely.  

Given a test subset $A$, let $N_A$ be the number of points in $A$ for
the random configuration (if the process is not simple, points are
counted with their multiplicities); this is a random variable, which
is defined on the state space and takes values in $\{0,1,\dots\}$. 

Assume that for any $A$, all the moments of $N_A$ are finite. Then we
assign to our point process a sequence $\rho_1,\rho_2,\dots$
of measures. The $n$th measure $\rho_n$ lives on the $n$-fold direct product
$X^n=X\times\dots\times  X$. In terms of the random configuration
$\xi=(x_1,x_2,\dots)$, $\rho_n$ is defined as the
expectation
$$
\rho_n={\Bbb E}\big\{
\sum_{i_1,\dots,i_n}\, \de(x_{i_1})\times\dots\times\de(x_{i_n})
\big\},
$$
where summation is taken over all $n$-tuples of 
{\it pairwise distinct\/} indices. Note that for any test subset $A$,
$$
\rho_n(A^n)={\Bbb E}\bigl\{N_A(N_A-1)\dots(N_A-n+1)\bigr\},
$$
the $n$th {\it factorial\/} moment of $N_A$. The measure $\rho_n$ is
called the $n$th {\it correlation measure\/} or the $n$th {\it
factorial moment measure\/}. 

Clearly, the correlation measures are symmetric and take finite
values on products of test subsets. The first correlation
measure $\rho_1$ is also called the {\it density measure\/}. The
value of $\rho_1$ on a test set $A$ is the mean number of points
(counted with multiplicities) occuring in $A$. 

Assume $X$ is equipped with a ``reference'' measure $dx$. When the
process is simple and the measure $\rho_n$ has a density 
$\rho_n(x_1,\dots,x_n)$ with respect to Lebesgue measure on $I^n$,
this density is called the $n$th {\it correlation function\/}.
Informally, $\rho_n(x_1,\dots,x_n)$ is equal to the probability that
the random configuration intersects each of infinitesimal volumes
$dx_1,\dots,dx_n$ around $x_1,\dots, x_n$, divided by 
$dx_1\dots dx_n$. When $X$ is a domain in an Euclidean space and $dx$
is Lebesgue measure it is convenient to regard the measure $\rho_n$
as a distribution and call the latter the correlation function, even
if we do not know {\it a priori\/} that the measure is absolutely
continuous with respect to $dx_1\dots dx_n$.

{}From now on, we take as $X$, the punctured interval
$I=[-1,1]\setminus\{0\}$. This is a locally compact space in the
natural topology, the point 0 playing the role of the infinity. A
subset of $I$ will be called ``bounded'' if it is relatively compact
in $I$, i.e., does not intersect a sufficiently small
interval $(-\ep,\ep)$. Thus, for a configuration in $I$, the only
possible accumulation point in ${\Bbb R}$ is 0.  

Define a map $\Om\to\Xi$ as 
$$
\om\;\mapsto\;\xi=(\al_1,\al_2,\dots,-\be_1,-\be_2,\dots),
$$
where all the $\al_i$ and $\be_i$ are assumed to be nonzero. In
particular, the point $\om=(\al,\be,\ga)=(0,0,1)$ is represented by
the empty configuration. One can verify that $\om\mapsto\xi$ is a
Borel map. Thus, any probability measure $P$ on $\Om$ defines a point
process in $I$ with state space $(\Om,\Xi)$; we shall denote this
process by $\PP$. 

\proclaim{Propostion 4.1} Any point process $\PP$ determined by a
probability measure $P$ on $\Om$ has the following special
property{\rm:} all the random variables of type $N_A$ are
bounded.
\endproclaim 

\demo{Proof} Take $\ep>0$ such that $A$ is contained in
$[-1,-\ep]\cup[\ep,1]$. Then, for any configuration $\xi$ originating
from a point $\om\in\Om$, 
$$
\sum_{x\in\xi\cap A} |x|\,\ge\,\ep\ |\xi\cap A|.
$$
Since the left--hand side does not exceed 1, we conclude that 
$|\xi\cap A|\le\ep^{-1}$. \qed
\enddemo

Consequently, $\PP$ possesses correlation measures.

\proclaim{Proposition 4.2} Consider a point process in $I$ possessing
correlation measures. If the diagonal in $I^2$ is a null set with
respect to $\rho_2$ then the process is simple. Conversely, if the
process is simple then, for any $n\ge2$ and any couple $i\ne j$ of
indices, the set
$$
\{(x_1,\dots,x_n)\in I^n\bigm| x_i=x_j\}
$$
is a null set with respect to $\rho_n$.
\endproclaim

\demo{Proof} This is verified by the same argument as in Proposition
3.5. Actually, the claim holds for general spaces $X$. \qed
\enddemo

\proclaim{Proposition 4.3} Let $P$ be a probability measure on $\Om$,
$\si_n$ its controlling measures, $\PP$ the point process defined by
$P$, and $\rho_n$ the correlation measures. Set
$$
(I^n)'=\{(x_1,\dots,x_n)\in I^n\bigm| 
x_i\ne x_j \quad\text{for $i\ne j$}\}.
$$

On $(I^n)'$, we have
$$
\rho_n=|x_1\dots x_n|^{-1}\si_n.
$$
\endproclaim

\demo{Proof} This is a direct consequence of the definitions of the
measures $\si_n$ and $\rho_n$. \qed
\enddemo

Thus, if $\PP$ is simple (which can be tested with the help of
Proposition 3.5), then the correlation measures are expressed in a
simple way through the controlling measures. Indeed, on the subset
$(I^n)'\subset I^n$ we can use the formula of Proposition 4.3, and the
complementary subset $I^n\setminus(I^n)'$ is negligible by virtue of
Proposition 4.2.

Let $\Pi(n,r)$ be the set of partitions of the {\it set\/}
$\{1,\dots,n\}$ consisting of $r$ nonempty blocs, $1\le r\le n$. An
element $\pi\in\Pi(n,r)$ can also be regarded as an equivalence
relation on $\{1,\dots,n\}$ with $r$ equivalence classes. We assign to
$\pi$ a subset $I^n_\pi\subset I^n$ obtained by intersecting $I^n$
with all hyperplanes of the form $x_i=x_j$ where $i\equiv j
\mod\pi$. Thus, $I^n_\pi$ is a `diagonal section' of $I^n$. Choosing a
representative $i$ in each equivalence class $\mod\pi$ we get natural
coordinates $\{x_i\}$ in the section $I^n_\pi$ by means of which we
can identify $I^n_\pi$ with $I^r$. Then we carry over the correlation
measure $\rho_r$ from $I^r$ to $I^n_\pi\subset I^n$ and denote the
resulting measure on $I^n$ by $\rho_{r,\pi}$.

\proclaim{Proposition 4.4} If $\PP$ is simple then
$$
\si_n\bigm|_{I^n}=|x_1\dots x_n|\,
\sum_{r=1}^n\, 
\sum_{\pi\in\Pi(n,r)}\rho_{r,\pi}\,.
$$
\endproclaim

\example{Example} Let $n=3$. There are 5 partitions $\pi$ of the set
$\{1,2,3\}$, and we have the following equality on $I^3$ (below 
$dx=dx_1dx_2dx_3$):
$$
\multline
\si_3(x_1,x_2,x_3)dx=|x_1x_2x_3|\{
\rho_3(x_1,x_2,x_3)dx\\
+\rho_2(x_1,x_2)\de(x_2-x_3)dx+
\rho_2(x_2,x_3)\de(x_3-x_1)dx+
\rho_2(x_3,x_1)\de(x_1-x_2)dx\\
+\rho_1(x_1)\de(x_1-x_2)\de(x_1-x_3)dx\}
\endmultline
$$
\endexample

\demo{Proof} Assume first that $P$ is the delta measure concentrated
at a point $\om=(\al,\be)\in\Om$ such that both $\al$ and $\be$ have
no repetitions. Let $\xi=(-\be_1<-\be_2<\dots<\al_2<\al_1)$ be the
corresponding subset of $I$. By the definition of the controlling
measures and the correlation measures,
$$
\gather
\si_n\bigm|_{I^n}=
\sum_{x_1,\dots,x_n\in\xi}
|x_1\dots x_n|
\de(x_1)\otimes\dots\otimes\de(x_n),\\
\rho_r=\sum
\Sb y_1,\dots,y_r\in\xi\\ y_i\ne y_j \endSb
\de(y_1)\otimes\dots\otimes\de(y_r).
\endgather
$$

To each $n$-tuple $(x_1,\dots,x_n)$ one can assign a partition $\pi$
of the set $\{1,\dots,n\}$ as follows: two distinct indices $i,j$ are
in the same bloc of $\pi$ if and only if $x_i=x_j$. Let $r$ be the
number of blocs of $\pi$. Then, in the sum over $(x_1,\dots,x_n)$, we
group together summands corresponding to the same bloc $\pi$. To pass
from $(x_1,\dots,x_n)$ to $(y_1,\dots,y_r)$ we order the blocs in an
arbitrary way and assign to the $i$th bloc $\pi_i$ the $i$th coordinate
$y_i$ (in other words, $y_i=x_j$ for any index $j\in\pi_i$). Then we
get the desired relation.

It is worth noting that in the above reasoning, enumeration of blocs
was used for convenience only: actually, the coordinates on the
section $I^n_\pi$ correspond just to blocs of $\pi$. 

Thus, we have verified the claim of the proposition in the particular
case when $P$ is a delta measure. In the general case, it remains  to
average over $\om$ with respect to $P$. \qed
\enddemo

{}From Proposition 3.5 or Proposition 4.4 it is clear that for
$n\ge2$, the controlling measure $\si_n$ can never be absolutely
continuous relative to Lebesgue measure on the cube $[-1,1]^n$,
because the diagonal sections $x_i=x_j$ of the cube always have a
nonzero mass. On the contrary, the correlation measures
$\rho_n$ can be absolutely continuous, at least, in the interior of
$I^n$. It will be shown in \cite{B} that the correlation functions
for $P=\Pz$ have analytic densities.  

\example{Remark 4.5} Even if a process of the form $\PP$ is not simple,
its correlation measures can be expressed through the controlling
measures. Conversely, each controlling measure $\si_n$ can be
expressed through the correlation measures $\rho_i$, $i=1,\dots,n$.  It
follows, in particular, that {\it the initial measure $P$ is uniquely
determined by the correlation measures of $\PP$\/}. Another proof of
this fact can be obtained from a general result of the theory of point
processes and Proposition 4.1.  \endexample
\medskip

Recall (see section 1) that the Thoma simplex possesses a natural
symmetry $\Om\to\Om$, 
$$
\om=(\al,\be)\,\mapsto\,\om^t=(\be,\al).\tag4.1
$$
Let $\PPz$ be the point process defined by the spectral measure $\Pz$
and $\rz_n$ be its correlation functions. 

\proclaim{Proposition 4.6} The symmetry map {\rm\tht{4.1}} takes
$\Pz$ to $P_{-z,-z'}$. Likewise, the symmetry $x\mapsto-x$ of $I$
takes the process $\PPz$ to the process $\Cal P_{-z,-z'}$. In
particular, we have 
$$
\gather
\si_n^{(-z,-z')}(x)=\siz_n(-x),\\
\rho_n^{(-z,-z')}(x)=\rz_n(-x).
\endgather
$$
\endproclaim

\demo{Proof} This follows from \tht{1.8}, \tht{2.2} and the
definition of the measures $\si_n$ and $\rho_n$. \qed
\enddemo

\head \S5. The density function \endhead

The aim of this section is to calculate the first correlation measure
(or the density measure) $\rz_1$ of $\PPz$, the point process
corresponding to $\Pz$. It will be shown that $\rz_1$ is absolutely
continuous with respect to Lebesgue measure on $I$. Hence, one can
speak about the first correlation {\it function\/} $\rz_1(x)$ which
is also called the density function. 

By Proposition 4.3, we have
$$
\rz_1=|x|^{-1}\siz_1\big|_I\,,
$$
so that it suffices to calculate the measure $\siz_1$ on $[-1,1]$.

\proclaim{Lemma 5.1} The measure $\siz_1$ is a unique solution of the
moment problem
$$
\int_{-1}^1x^l\siz_1(dx)=
\sum\Sb p,q\ge0\\ p+q=l\endSb
\frac{(-1)^q\,t\,(z+1)_p(-z+1)_q(z'+1)_p(-z'+1)_q}
{(t)_{p+q+1}\,(p+q+1)\,p!\,q!}\,,
\tag5.1
$$
where $l=0,1,2,\dots$.
\endproclaim

\demo{Proof} We apply Proposition 3.3. Since $n=1$, the parameter $d$
takes the only value 1, and the summation is taken over the hook
diagrams $\la=(p\,|\,q)$.

It follows from the Murnaghan--Nakayama rule \cite{M, \S{I}.7, Ex. 5}
that
$$
\chi^{(p\,|\,q)}_{(p+q+1)}=(-1)^q.
$$
This yields the desired formula. As was already mentioned, uniqueness
holds because the support of the measure in question is bounded. \qed
\enddemo

We shall deal with the distributions
$$
\phi_a(u)=\frac{u_+^a}{\Ga(a+1)}\,,\qquad u\in\R,\quad a\in\C,
$$
concentrated on the right semiaxis $\R_+$. Here the numerator $u_+^a$
coincides with the function $u^a$ when $u>0$ and vanishes when
$u<0$. When $\Re a>-1$, $\phi_a$ is an integrable function, and when
$\Re a\le-1$, it is defined via analytic continuation. \footnote{Here
and in what follows, to simplify the notation, we shall write
distributions as if they were ordinary functions.} For any 
$a\in\C$, the only singularity of the distribution $\phi_a$ may be at
0. 

The product
$$
\phi_{ab}(u)=\phi_a(u)\,\phi_b(1-u),\qquad a,b\in\C
$$
is correctly defined because the possible singularities of the
factors are at different points ($u=0$ and $u=1$,
respectively). The result is a distribution concentrated on $[-1,1]$.

The following formula holds:
$$
\int u^p(1-u)^q\phi_{ab}(u)=
\frac{(a+1)_p(b+1)_q}{\Ga(a+b+p+q+2)}\,, \qquad p,q=0,1,2,\dots
\tag5.2
$$
Indeed, when $\Re a>-p-1$ and $\Re b>-q-1$, this formula is
equivalent to the classical Euler beta integral formula, and for
arbitrary $a,b\in\C$ the result holds by analytic continuation.

Set
$$
\Phi(a+1;\ze)=\int e^{\ze u} \phi_{a,-a}(u)du,\qquad \ze\in\C.
$$
I.e., $\Phi(a+1;\,\cdot\,)$ is the Laplace transform of
$\phi_{a,-a}$. Expanding the exponential function and using the beta
integral \tht{5.2} with $q=0$, one sees that 
this is a special case of Kummer's hypergeometric function $_1 F_1$
(also called the confluent hypergeometric function):
$$
\Phi(a+1;\ze)=\sum_{k=0}^\infty
\frac{(a+1)_k}{(k+1)!\,k!}\,\ze^k={_1F_1}(a+1;2;\ze),
$$
see \cite{E, ch. 6}. Clearly, $\Phi(a+1;\ze)$ is an entire function in
$\ze$.

\proclaim{Theorem 5.2} The density measure satisfies the following
equation 
$$
\gathered
\int_{-1}^1\Phi(t+1;\ze x)\,\siz_1(dx)=
\Phi(z+1;\ze)\,\Phi(-z'+1;-\ze)\\
=\Phi(z'+1;\ze)\,\Phi(-z+1;-\ze), \qquad \ze\in\C.
\endgathered
\tag5.3
$$
\endproclaim

\demo{Proof} We rewrite formula \tht{5.1} as
$$
\gather
\frac{\ze^l(t+1)_l}{(l+1)l!}\,
\int_{-1}^1x^l\siz_1(dx)\\
=\sum\Sb p,q\ge0\\ p+q=l\endSb
\frac{(z+1)_p(-z+1)_q}{(p+q+1)!}\,
\frac{(z'+1)_p(-z'+1)_q}{(p+q+1)!}\,
\frac{\ze^p(-\ze)^q}{p!\,q!}\,.
\endgather
$$
Next, we replace the first and the second ratios on the right by the
corresponding beta integrals \tht{5.2} and sum over $l$. The result
will be as follows
$$
\gather
\sum_{l\ge0}\int_0^1\frac{(\ze x)^l(t+1)_l}{(l+1)!\,l!}\,
\siz_1(dx)\\
=\sum_{p,q\ge0}\iint
\phi_{z,-z}(u)\phi_{z',-z'}(v)
\frac{(\ze uv)^p}{p!}
\frac{(-\ze(1-u)(1-v))^q}{q!}\,dudv.
\endgather
$$
Interchanging summation and integration, we get
$$
\int_0^1\Phi(t+1;\ze x)\,\siz_1(dx)=
\iint e^{\ze(uv-(1-u)(1-v))}\phi_{z,-z}(u)\phi_{z',-z'}(v)dudv.
$$

Since
$$
uv-(1-u)(1-v)=u+v-1=u-(1-v)=v-(1-u),
$$
the double integral on the right factorizes into a product of two
one--dimensional integrals. Making a change of a variable,
$(1-v)\mapsto v$ or $(1-u)\mapsto u$, and using the definiton of the
function $\Phi$, we obtain the first and the second variants of
formula \tht{5.3}, respectively.

Note that the equivalence of both variants of formula \tht{5.3} also
follows from the identity
$$
\Phi(a+1;\ze)=e^\ze\Phi(-a+1;\ze),
$$
which is a particular case of Kummer's transform for $_1F_1$, see
\cite{E, \S6.3, (21)}. \qed
\enddemo

\example{Remark 5.3} The idea to use the beta integral is due to
Borodin. My initial proof of Theorem 5.2 was more complicated: I
transformed the right--hand side of \tht{5.1} to the form
$$
\frac{l!\,(l+1)!\,t}{(t)_{l+1}}\,
\sum_{p+q=l}\frac{(-z+1)_q(z'+1)_p}{p!\,(p+1)!\,q!\,(q+1)}
$$
with the help of a Leibniz--type formula for the difference operator
$f(z)\mapsto f(z)-f(z-1)$. \qed
\endexample
\medskip

Set
$$
\phi_{ab}^{(-)}(u)=\phi_{ab}(-u);
$$
this is a distribution concentrated on $[-1,0]$.

\proclaim{Corollary 5.4} If $t=1$ then the measure $\siz_1$ is the
convolution product of two distributions, concentrated on $[0,1]$ and
$[-1,0]${\rm:} 
$$
\siz_1=\phi_{z,-z}\,*\,\phi_{-z',z'}^{(-)}.
$$
\endproclaim

\demo{Proof} When $t=1$, the function $\Phi(t+1;\ze)$ degenerates to
the exponential $e^\ze$. It follows that the left--hand side of
\tht{5.3} reduces to the Laplace transform of the measure $\siz_1$.
On the other hand, the right--hand side of \tht{5.3} is the product
of the Laplace transforms of the distributions $\phi_{z,-z}$ and
$\phi_{-z',z'}^{(-)}$. \qed
\enddemo

In the general case, to extract from the equation \tht{5.3} an
expresssion for $\siz_1$ we need a little formalism which will also be
employed in \S6.

Define an operation $\odot$ on distributions by the rule
$$
(A\odot B)(x)=\iint\de(x-y_1y_2)A(y_1)B(y_2)dy_1dy_2.
$$
Or, in terms of a test function $\psi$,
$$
\int(A\odot B)(x)\psi(x)dx=\iint A(y_1)B(y_2)\psi(y_1y_2)dy_1dy_2.
$$
The operation $\odot$ may be viewed as the convolution product on the
semigroup $(\R,\,\cdot\,)$ of real numbers under multiplication; for
this reason we shall call it {\it pseudoconvolution\/}. 

Of course, in order that the pseudoconvolution $A\odot B$ be
correctly defined, $A$ and $B$ must satisfy appropriate conditions.
For example, it suffices that they would be compactly supported; then
$A\odot B$ will be compactly supported, too; moreover, 
$\operatorname(A\odot B)\subseteq\operatorname{supp}
A\,\cdot\,\operatorname{supp}B$. We shall employ this operation
for distributions concentrated on $[-1,1]$; since
$[-1,1]$ is a subsemigroup of $(\R,\,\cdot\,)$, the result will
always be a distribution of the same kind.

\example{Remark 5.5} (i) Assume that $A$ and $B$ are integrable
functions or (complex) measures of finite variance with no atom at 0.
Then the point 0 may be neglected, and the pseudoconvolution reduces,
in essence, to the conventional convolution on the multiplicative
group $\R^*=\R\setminus\{0\}$.

(ii) The same is also true if both $A$ and $B$ are distributions
whose supports do not contain 0. However, in the general case, the
point 0 can cause complications.

(iii) Here is an illustrative example of what can happen in the  
extreme case when both $A$ and $B$ are supported at 0: Denoting by
$\de_0$ the delta function and by $\de_0^{(m)}$ its derivative of
order $m$, we have
$$
\de_0^{(m)}\odot\de_0^{(n)}=\cases 
m!\,\de_0^{(m)}, &m=n\\ 0,&m\ne n. \endcases
$$
\endexample
\medskip

\proclaim{Lemma 5.6} Formula \tht{5.3} of Theorem 5.2 is equivalent
to 
$$
\siz_1\odot\phi_{t,-t}=\phi_{z,-z}*\phi_{-z',z,}^{(-)}.
\tag5.4
$$
\endproclaim

\demo{Proof} We claim that \tht{5.3} coincides with the Laplace
transform of the latter formula. Indeed, for the right--hand side
this was already noted in the proof of Corollary 5.3, and for the
left--hand side, this is readily verified by substituting the
integral reprsentation of the function $\Phi(t+1;\,\cdot\,)$ to the
left--hand side of \tht{5.3}. \qed
\enddemo

\proclaim{Lemma 5.7} Denoting by $\de_1$ the Dirac mass at the point
1, we have
$$
\phi_{t,-t}\odot\phi_{1,t-2}=\dfrac 1{\Ga(t+1)}\de_1.
$$
\endproclaim

\demo{Proof} Let us check that the $n$th moment of left--hand side
(where $n=0,1,2,\dots$) is equal to $(\Ga(t+1))^{-1}$. By the
definition of $A\odot B$, the $n$th moment of $A\odot B$ is equal to
the product of the $n$th moments of $A$ and $B$. Next, observe that
$$
\int x^n\phi_{ab}(x)dx=\frac{(a+1)_n}{\Ga(a+b+n+2)}\,.
$$
Therefore, the $n$th moment in question is equal to
$$
\frac{(t+1)_n}{\Ga(n+2)}\,\frac{(2)_n}{\Ga(t+n+1)}
=\frac 1{\Ga(t+1)}\,,
$$
as was to be shown. \qed
\enddemo

\proclaim{Theorem 5.8} We have
$$
\siz_1=\Ga(t+1)\,(\phi_{z,-z}*\phi_{-z',z'}^{(-)})\odot\phi_{1,t-2}\,,
\tag5.5
$$
so that
$$
\rz_1(x)=\frac{\Ga(t+1)}{|x|}\,
\{(\phi_{z,-z}*\phi_{-z',z'}^{(-)})\odot\phi_{1,t-2}\}\,,
\qquad x\ne0.
\tag5.6
$$
\endproclaim

\demo{Proof} Take the pseudoconvolution of the both sides of
formula \tht{5.4} of Lemma 5.6 with $\phi_{1,t-2}$ and apply then
Lemma 5.7. \qed
\enddemo

The expression of Theorem 5.8 is a two--dimensional integral
representation of the measure $\siz_1$. It can also be derived from
the moment formula \tht{5.1} by the same method as that we employed in
the proof of Theorem 5.2.

\example{Remark 5.9} There is a somewhat different formula for
$\siz_1$: set $D=x\frac d{dx}$; then
$$
\siz_1(x)=\Ga(t+1)\,(-D)\,
[(\phi_{z,-z}*\phi_{-z',z,}^{(-)})\odot\phi_{0,t-1} ].
$$
Indeed, by virtue of Proposition 6.8 (see below), the right--hand
side is equal to 
$$
\Ga(t+1)\,(\phi_{z,-z}*\phi_{-z',z'}^{(-)})
\odot(-D\phi_{0,t-1}).
$$ 
It is readily verified that 
$$
-D\,\phi_{0,t-1}=\phi_{1,t-1}.
$$ 
So, the formula in question is equivalent to \tht{5.5}.
\endexample
\medskip

\proclaim{Theorem 5.10} The restriction of the density measure
$\rz_1$ to $(-1,0)\cup(0,1)$ is absolutely continuous with respect to
Lebesgue measure $dx${\rm:}
$$
\rz_1(dx)=\rz_1(x)dx,
$$
and the density function $\rz_1(x)$ is real analytic on $(0,1)$ and on
$(-1,0)$.  

When $x\in(0,1)$, the density function can be written in the
following form
$$
\rz_1(x)=\frac{\Ga(t+1)}{\Ga(z+1)\Ga(z'+1)}\,
(1-x)^{-z-z'+t}\,
\langle A,\Psi\rangle,
\qquad x\in(0,1),
\tag5.7
$$
where $A$ is a two--dimensional distribution concentrated on the
triangle $u,v\ge0$, $u+v\le1$, and $\Psi$ is a test function, which
is correctly defined and smooth in a neighborhood of that
triangle{\rm:} 
$$
\gather
A=A(u,v)=\phi_{-z}(u)\phi_{-z'}(v)\phi_{t-2}(1-u-v),\\
\Psi=\Psi(u,v)=
(1-(1-x)u)^z(1-(1-x)v)^{z'}(1-(1-x)(u+v))^{-t}.
\endgather
$$

When $x\in(-1,0)$, the density function can be written in the same
form, with the only modification: $x$ is replaced by $|x|$ and the
parameters $z,z'$ are replaced by $-z,-z'$.
\endproclaim

\demo{Proof} We start with a formal transformation of the expression
\tht{5.6} given by Theorem 5.8. 

Assume that $x\in(0,1)$. By \tht{5.6} and the definition of $\odot$,
$$
\siz_1(x)=\Ga(t+1)\,\iiint\de(x-w(u+v-1))
\phi_{z,-z}(u)\phi_{z',-z'}(v)\phi_{1,t-2}(w)dudvdw,
$$
where $\de(\,\cdot\,)$ is the delta function. We can exclude the
integration over $w$ using the formula 
$$
\gather
\frac 1x \,\de(x-w(u+v-1))\phi_{1,t-2}(w)dw\\
=\frac 1{x(u+v-1)}\,\de(w-\frac x{u+v-1})\,\phi_{1,t-2}(w)dw\\
=\frac 1{x(u+v-1)}\,\phi_{1,t-2}(\frac x{u+v-1})\\
=\frac 1{(u+v-1)^t}\,\phi_{t-2}(u+v-1-x).
\endgather
$$
This gives
$$
\gather
\rz_1(x)=\Ga(t+1)\iint 
\phi_z(u)\phi_{-z}(1-u)\phi_{z'}(v)\phi_{-z'}(1-v)\\
\times\phi_{t-2}(u+v-1-x)\frac{dudv}{(u+v-1)^t}\,.
\endgather
$$
Make a change of variables,
$$
1-u=(1-x)u_1,\quad 1-v=(1-x)v_1.
$$
Then
$$
\gather
dudv=(1-x)^2du_1dv_1,\\
u=1-(1-x)u_1,\quad v=1-(1-x)v_1,\\
u+v-1=1-(1-x)(u_1+v_1),\quad u+v-1-x=(1-x)(1-u_1-v_1).
\endgather
$$
Substituting these expressions into the integral and renaming then
the variables $u_1,v_1$ to $u,v$, we come to formula \tht{5.7}.
In the case $x\in(-1,0)$ one can use exactly the same argument;
an alternative possibility is to use the symmetry property indicated
in Lemma 4.6.

To justify these formal transformations from \tht{5.6} to \tht{5.7},
let us, for a moment, interpret $z$, $-z$, $z'$, $-z'$, and $t-2$ as
five independent complex parameters. If the real parts of these
variables are strictly positive, then all the distributions become
ordinary continuous functions, and our transformations are readily
justified. On the other hand, the expressions \tht{5.6} (with $x>0$)
and \tht{5.7} are both correctly defined distributions in $x$, which
depend holomorphically on the parameters. Hence, by the principle
of analytic continuation, they are equivalent.

Finally, the expresion \tht{5.7} is real analytic, because the test
function $\Psi$ is correctly defined for complex values of the
parameters $x$ in the strip $0<\Re x<1$ and depends on $x$
analytically. \qed
\enddemo

\example{Remark 5.11} By Lemma 6.10 (see below), the density
measure $\rz_1$ has no atom at 0, and one can show that $\rz_1$ has no
atoms at the end points $\pm1$. So, it is completely determined by
the density function on the open intervals $(0,1)$ and $(-1,0)$.
\endexample
\medskip

We shall express the density function through a multivariate
hypergeometric series. Let 
$$
a=(a_1,\dots,a_n),\quad b=(b_1,\dots,b_n),\quad c
$$
be complex parameters and 
$$
y=(y_1,\dots,y_n)
$$
be $n$ complex variables. The $n$-dimensional {\it Lauricella
hypergeometric function\/} of type $B$ is defined by the series
$$
F_B^{[n]}(a,b\,;c\,|\,y)=
\sum_{m_1,\dots,m_n\ge0}
\frac{(a_1)_{m_1}(b_1)_{m_1}\dots (a_n)_{m_n}(b_n)_{m_n}}
{(c)_{m_1+\dots m_n}\,m_1!\dots m_n!}\,
y_1^{m_1}\dots y_n^{m_n}\,
$$
where the series is absolutely convergent for $|y_1|<1, \dots,
|y_n|<1$, see \cite {AK}, \cite{Ex}. When
$n=1$, this is Gauss' hypergeometric function, and when $n=2$, this
is Appell's hypergeometric function $F_3$. Note also that the
function $F_B$ remains invariant when the couples
$(a_1,b_1),\dots, (a_n,b_n)$ are permuted or the parameters in
couples are interchanged. 

\proclaim{Theorem 5.12} Let
$$
\gather
a=(-z,\,-z',\,t-1),\quad b=(-z+1,\,-z'+1,\,t),\\
c=-z-z'+t+1=(-z+1)(-z'+1).
\endgather
$$

For $0<x\le1$, we have
$$
\rz_1(x)=\frac{\Ga(t+1)}{\Ga(z+1)\Ga(z'+1)}\,
\phi_{c-1}(1-x)\, F_B^{[3]}(a,b\,;c\,|\,1-x,1-x,1-x). \tag5.8
$$
The same expression holds for $-1\le x<0$ provided that $x$ is
replaced by $|x|$ and $z,z'$ are multiplied by $-1$.
\endproclaim

\demo{Proof} First of all, note that the parameter $c$, as defined
above, is a strictly positive real number (this follows from the
fundamental assumptions on the parameters $z,z'$). It follows that
$\phi_{c-1}(1-x)$ is integrable at $x=1$. Consequently, the whole
expression is integrable at $x=1$, as it should be. 

As before, the case of negative $x$ can be reduced to that of positive
$x$ by symmetry. Next, by Remark 5.11, we may suppose $0<x<1$.

Set
$$
\gather
u_1=u,\quad u_2=v,\quad u_3=u+v,\\
d_1=-z,\quad d_2=-z',\quad d_3=t-2,\\
y_1=y_2=y_3=1-x.
\endgather
$$
By the binomial expansion,
$$
\Psi=\sum_{m_1,m_2,m_3\ge0}
\frac{(-z)_{m_1}(-z')_{m_2}(t)_{m_3}}
{m_1!\,m_2!\,m_3!}\,
(y_1u_1)^{m_1}(y_2u_2)^{m_2}(y_3u_3)^{m_3}
$$
Substituting this in \tht{5.7} and
employing the Euler beta integral
$$
\gather
\int\limits\Sb u_1+u_2+u_3=1\\ u_1,u_2,u_3\ge0\endSb
\phi_{d_1}(u_1)\,\phi_{d_2}(u_2)\,\phi_{d_3}(u_3)\,
u_1^{m_1}\,u_2^{m_2}\,u_3^{m_3}\,du\\
=\frac{(d_1+1)_{m_1}(d_2+1)_{m_2}(d_3+1)_{m_3}}
{\Ga(d_1+d_2+d_3+3)\,(d_1+d_2+d_3+3)_{m_1+m_2+m_3}}\,
\endgather
$$
we obtain the desired formula. \qed
\enddemo

\head \S6. An application  \endhead

In this section we shall prove the following result. \footnote{I am
grateful to Jean--Louis Clerc and Jacques Faraut for discussions
related to the proof of Theorem 6.1.}

\proclaim{Theorem 6.1} All the spectral measures $\Pz$ are
concentrated on the face $\Om_0$ defined in \tht{1.5}. 
\endproclaim

By Proposition 3.4, this is equivalent to the fact that $\siz_1$ has
no atom at 0. The proof of the latter claim is divided into a series of
lemmas. 

Let $C_0(\R)$ be the space of continuous compactly supported
functions on $\R$ and $C_0^\infty (\R)$ be its subspace consisting of
smooth compactly supported functions. Let $a,b$ be complex parameters
and $x$ be the coordinate on $\R$. Let 
$$
\Fp_a=\phi_a\cdot C_0^\infty(\R)
$$
denote the space of distributions formed by the products
$\phi_a\cdot f$ with $f\in C_0^\infty(\R)$, and let $\Fm_a$ denote
the image of $\Fp_a$ under the reflection $x\mapsto -x$.

\proclaim{Lemma 6.2} We have
$$
\frac d{dx}\,\Fpm_a\subset\Fpm_{a-1}+\Fpm_a.
$$
\endproclaim

\demo{Proof} Indeed, this follows from the well--known formula
$$
\frac d{dx}\,\phi_a(x)=\phi_{a-1}(x).
$$
\qed
\enddemo

\proclaim{Lemma 6.3} We have
$$
\Fpm_a\subset \frac d{dx}\,\Fpm_{a+1}+\Fpm_{a+1}.
$$
\endproclaim

\demo{Proof} Indeed, for any $f\in C_0^\infty(\R)$,
$$
\phi_a f=(\phi_{a+1}f)'-\phi_{a+1}f'
\in\,(\frac d{dx}\,\Fp_{a+1}\,+\,\Fp_{a+1}).
$$
The same argument works for the sign ``$-$''.  \qed
\enddemo

\proclaim{Lemma 6.4} We have
$$
\Fp_a\,*\,\Fm_b\subset
(\Fp_{a-1}\,*\,\Fm_{b+1})
\,+\,
(\Fp_a\,*\,\Fm_{b+1}).
$$
\endproclaim

\demo{Proof} By Lemma 6.3,
$$
\gather
\Fp_a\,*\,\Fm_b\subset
\Fp_a\,*\,(\frac d{dx}\,\Fm_{b+1}\,+\,\Fm_{b+1})\\
=(\Fp_a\,*\,\frac d{dx}\,\Fm_{b+1})\,+\,
(\Fp_a\,*\,\Fm_{b+1})\\
=(\frac d{dx}\,\Fp_a\,*\,\Fm_{b+1})\,+\,
(\Fp_a\,*\,\Fm_{b+1}).
\endgather
$$
Next, we apply Lemma 6.2. \qed
\enddemo

\proclaim{Lemma 6.5} Assume that one of the following conditions
holds: 

{\rm (i)} $\Re a>0$, $\Re b<1$;

{\rm (ii)} $\Re a=\Re b=0$.

Then $\Fp_a\,*\,\Fm_{-b}\subset C_0(\R)$.
\endproclaim

\demo{Proof} (i) In this case, the elements of $\Fp_a$ are continuous
functions with compact support while the elements of $\Fm_{-b}$ are
integrable functions with compact support. Therefore, the result of
convolution are continuous functions. 

(ii) In this case, the elements of the both spaces are bounded
measurable functions with compact support. So, they are
square integrable functions. Therefore, the result of convolution are
again continuous functions. \qed
\enddemo

\proclaim{Lemma 6.6} Assume that the parameters $a,b$ satisfy the
following condition{\rm:} there exist integers $m\ge n\ge0$ such that
either 
$$
m<\Re a<m+1,\qquad n<\Re b<n+1
$$
or
$$
\Re a=m, \qquad \Re b=n.
$$
Then 
$$
\Fp_a\,*\,\Fm_{-b}\subset C_0(\R).
$$
\endproclaim

Actually, we need only the case $m=n$, but it will be convenient to
check a slightly more general claim with $m\ge n$. 

\demo{Proof} Assume first that $n\ge1$. By lemma 6.4,
$$
\Fp_a\,*\,\Fm_{-b}\subset
(\Fp_{a-1}\,*\,\Fm_{-(b-1)})\,+\,
(\Fp_a\,*\,\Fm_{-(b-1)}),
$$
which enables one to reduce the claim of the lemma for a given
couple $m\ge n\ge1$ to the same claims with $(m,n)$ replaced by
$(m-1,n-1)$ or by $(m,n-1)$. The reduction stops when $n=0$, but then
we can apply Lemma 6.5. \qed
\enddemo 

\proclaim{Lemma 6.7} Near $0$, the distribution
$$
\phi_{z,-z}\,*\,\phi_{-z',z'}^{(-)}
$$ 
is given by a continuous function.
\endproclaim

\demo{Proof} By making use of an appropriate partition of unity one
can represent the both factors in the form 
$$
\phi_{z,-z}=A_0+A+A_1,\qquad
\phi_{-z',z'}^{(-)}=B_{-1}+B+B_0,
$$
where $A,B$ are smooth functions and $A_0,A_1,B_{-1}, B_0$ are
certain distributions concentrated near the points 0, 1, $-1$, 0,
respectively (they are obtained by multiplying the initial
distributions $\phi_{z,-z}$ and $\phi_{-z',z'}^{(-)}$ by appropriate
functions from $C_0^\infty(\R)$). It follows that 
$$
\phi_{z,-z}\,*\,\phi_{-z',z'}^{(-)}
=(A_0*B_0)+(A_1*B_{-1})+(A_0*B_{-1})+(A_1*B_0)+(\dots),
$$
where $(\dots)$ is a smooth function. Since the distributions
$A_0*B_{-1}$ and $A_1*B_0$ are concentrated near the points $-1$ and
1, respectively, it remains to check that the distributions
$A_0*B_0$ and $A_1*B_{-1}$ are actually continuous functions.

We have $A_0\in\Fp_z$ and $B_0\in\Fm_{-z'}$, so that
$$
A_0\,*\,B_0\,\in\,(\Fp_z\,*\,\Fm_{-z'}).
$$
Next, the distrubution $A_1*B_{-1}$ will not change if we shift the
both factors by 1, to the left and to the right, respectively. The
resulting distributions will lie in $\Fm_{-z}$ and $\Fp_{z'}$,
respectively, whence
$$
A_1\,*\,B_{-1}\,\in\,(\Fp_{z'}\,*\,\Fm_{-z}).
$$

Now, let us compare our fundamental assumptions on the parameters
$z,z'$ with the assumptions on the parameters $a,b$ in Lemma 6.6.
Without loss of generality, we may assume that $\Re z\ge0$ and $\Re
z'\ge0$ (otherwise we may apply Proposition 4.6). Then the couples
$(a,b)=(z,z')$ and $(a,b)=(z',z)$ will satisfy the assumptions of
Lemma 6.6 with $m=n$. Application of this lemma concludes the proof.
\qed 
\enddemo

\proclaim{Lemma 6.8} Set $D=x\dfrac d{dx}$.
For any compactly supported distributions $A,B$, we have
$$
(DA)\odot B=D(A\odot B)=A\odot(DB).
$$
\endproclaim

Note that the claim is obvious when the supports of $A$ and $B$ do not
contain 0, because then, replacing $x$ by $s=\log|x|$, we can reduce
the pseudoconvolution to the ordinary convolution and the operator $D$
--- to $\dfrac d{ds}$.

\demo{Proof} Set $D'=\dfrac d{dx}\circ x=D+1$. For an arbitrary
smooth test function $\psi$, we have
$$
\gather
\langle(DA)\odot B,\psi\rangle
=\iint DA(y_1)B(y_2)\psi(y_1y_2)dy_1dy_2\\
=\int B(y_2)dy_2\int DA(y_1)\psi(y_1y_2)dy_1\\
=\int B(y_2)dy_2\int A(y_1)(-y_1\dfrac d{dy_1}-1)\psi(y_1y_2)dy_1\\
=\int B(y_2)dy_2\int A(y_1)(-D'\psi)(y_1y_2)dy_1\\
=\langle A\odot B,-D'\psi\rangle =
\langle D(A\odot B),\psi\rangle.
\endgather
$$
This proves the first equality, and the second one is verified
similarly. \qed
\enddemo

\proclaim{Lemma 6.9} Assume that $A$ is a distribution with a compact
support not containing 0, and $B$ belongs to a class $\Fp_a$ with 
$\Re a>-1$. Then the distribution $A\odot B$ is a function. 
\endproclaim

\demo{Proof} Since the support of $A$ is separated from 0, $A$ can be
represented in the form $A=f(D)\widetilde A$, where $f$ is a
polynomial, $D$ is as in Lemma 6.8, and $\widetilde A$ is an ordinary
(say, integrable) function. By Lemma 6.8,
$$
A\odot B=(f(D)\widetilde A)\odot B=\widetilde A\odot(f(D)B).
$$
On the other hand, we have
$$
D\phi_a=a\phi_a,\qquad a\in\C,
$$
which implies
$$
D\Fp_a\subseteq\Fp_a,\qquad a\in\C,
$$
so that $f(D)B$ lies in $\Fp_a$ together with $B$. When $\Re a>-1$,
the distributions of class $\Fp_a$ are ordinary (integrable)
functions. Consequently, the pseudoconvolution of $\widetilde A$ and
$f(D)B$ is, in essence, the ordinary convolution product on the
multiplicative group $\R^*$ (see Remark 5.5), and its result is an
ordinary function. \qed
\enddemo

The next lemma concludes the proof of Theorem 6.1.

\proclaim{Lemma 6.10} The measure $\siz_1$ has no atom at $0$. 
\endproclaim

\demo{Proof} We shall use different arguments when $t<1$, $t=1$, and
$t>1$. 

When $t=1$, the measure $\siz_1$ is given by the expression 
of Corollary 5.4, and the claim follows from Lemma 6.7.

Assume that $t<1$ and examine the expression \tht{5.4}. Since
$0<t<1$, the distribution $\phi_{t,-t}$ is actually a nonnegative
integrable function. We can decompose the measure $\siz_1$ into the
sum of two components: one is $\operatorname{const}\cdot\de_0$ (a
multiple of the Dirac mass at 0) and another is a measure $A$ with no
atom at 0. The measure $A\odot\phi_{t,-t}$ cannot have an atom at 0
(see Remark 5.5). On the other hand, 
$$
\operatorname{const}\de_0\odot \phi_{t,-t}
=\operatorname{const}\,(\int\phi_{t,-t}(y)dy)\,\de_0,
$$
where the integral is strictly positive. But the right--hand side of
formula \tht{5.4} cannot have an atom at 0, by virtue of Lemma 6.7.
Hence, $\operatorname{const}=0$, so that $\siz_1$ has no atom at 0,
too.

Finally, assume that $t>1$ and look at formula \tht{5.5}. By Lemma
6.7, one can write
$$
\Ga(t+1)\,\phi_{z,-z}*\phi_{-z',z'}^{(-)}=A_0+A_1,
$$
where $A_0$ is a continuous function and $A_1$ is a distribution whose
support does not contain 0. Next, since $t>1$, the distribution
$\phi_{1,t-2}$ can be written as 
$$
\phi_{1,t-2}=B_0+B_1,
$$
where $B_0\in\Fp_1$ and $B_1$ is an integrable function concentrated
near 1. Then we have
$$
\gather
\siz_1=(A_0+A_1)\odot(B_0+B_1)\\
=(A_0\odot B_0)+(A_1\odot B_0)+(A_0\odot B_1)+(A_1\odot B_1).
\endgather
$$

Let examine the four summands of the latter expression. Since $A_0$,
$B_0$, $B_1$ are ordinary (integrable) functions, the terms 
$A_0\odot B_0$ and $A_0\odot B_1$ are ordinary functions according to
Remark 5.5. The term $A_1\odot B_0$ is an ordinary function by Lemma
6.9. The term $A_1\odot B_1$ is a distribution concentrated outside a
neighborhood of 0. We conclude that the whole expression cannot have an
atom at 0. \qed 
\enddemo

\head 7. Appendix: Correlation functions of Poisson--Dirichlet 
processes \endhead

In this appendix, we briefly discuss a parallel (but more simple)
theory. It is related to the Poisson--Dirichlet distributions, a
remarkable one--parametric family of probability measures
$\{PD(t)\,|\,t>0\}$ which live on an infinite--dimensional subsimplex
$\Delta\subset \Omega$. The measures $PD(t)$ determine point processes
$\Cal{PD}(t)$ on $(0,1]$, called the Poisson--Dirichlet processes, and
we shall calculate the correlation functions of $\Cal{PD}(t)$ by using
the general formalism of \S\S1--4.

The Poisson--Dirichlet distributions were studied in many papers from 
different points of view, see, e.g., \cite{Ki2, W}. Our interest in them 
is caused by the fact that they are one of the basic elements in the 
construction \cite{KOV} of the generalized regular representations $T_z$. 
So, both kind of measures, the $PD(t)$'s and the $\Pz$'s, are connected 
with the same construction  --- that of the representations $T_z$. 
But they appear at different levels of that construction: the 
former --- at the `group' level (as they are responsible for certain 
quasiinvariant measures on  $\widetilde{G/K}$, see the Introduction), 
and the later --- at the `dual' level (the representation level). 
However, both $PD(t)$ and $\Pz$  are spectral measures in the sense 
that they govern decomposition of certain objects into indecomposable 
ones\footnote{At the `group' level, the `objects' are $K$-invariant 
probability measures on $\widetilde{G/K}$ and the `indecomposable 
objects' are ergodic ones, see \cite{Ki1, KOV}.}, and both $PD(t)$ 
and $\Pz$ can be interpreted as point processes.

For these reasons, it seems interesting to compare the point processes 
$\Cal{PD}(t)$ and $\PPz$, and the main purpose of the present 
appendix is to prepare a foundation for such a comparison (we postpone 
the discussion to the third article \cite{BO}). Another purpose is to 
illustrate the formalism of sections 1--4 on a simpler material.

We define the simplex $\De$ as the closed subset of $\Omega$ determined 
by $\beta_1=\beta_2=\dots=0$. In other words, $\De$ is the set of 
sequences $\al=(\al_1\ge\al_2\ge\dots\ge0)$ such that $\sum\al_i\le1$. 
The specialization \tht{1.3} is replaced by the following one:
$$
\gather
p_1\;\mapsto\;\wtp_1=\wtp_1(\al;\ga):=\sum_{i=1}^\infty\al_i+\ga\\
p_n\;\mapsto\;\wtp_n=\wtp_n(\al;\ga):=\sum_{i=1}^\infty\al_i^n,
\qquad n\ge0,
\endgather
$$
which is equivalent to
$$
1+\sum_{n=1}^\infty h_nu^n\;\mapsto\;
1+\sum_{n=1}^\infty \wth_n u^n:=e^{\ga u}
\prod_{i=1}^\infty \frac 1{1-\al_i u}\,.
$$

In what follows, we shall assume that $\al$ is a point of $\Delta$ and 
$$
\ga=1-\sum_{i=1}^\infty\al_i\,.
$$
Then all $\wtp_n$ turn into continuous functions functions on $\De$ 
(note that $\wtp_1\equiv1$). Consequently, any element $f\in\La$ is 
converted into a continuous function $\wtf=\wtf(\al)$ on $\De$.

In place of the Schur functions $s_\mu$ we shall deal with the 
{\it monomial symmetric functions\/} $m_\la\in\La$. Recall that in 
the standard realization of $\La$ as the algebra of symmetric 
functions in variables $x_1,x_2,\dots$, the function $m_\la$ is the 
sum of all {\it distinct\/} monomials obtained from
$x^\la=x_1^{\la_1}x_2^{\la_2}\dots$ by permuting variables, 
see \cite{M}.

The Pieri formula \tht{1.3} for the Schur functions is replaced by 
its counterpart for the monomial symmetric functions,
$$
m_\mu\cdot p_1=\sum_{\la:\,\la\searrow\mu}\ka_0(\mu,\la)\,m_\la,
$$
where the coefficients $\ka_0(\mu,\la)$ are positive integers defined 
as follows. Given $\la\searrow\mu$, there exists a unique 
$i\in\{1,2,\dots\}$ such that $\la_i=\mu_i+1$ (and $\la_j=\mu_j$ 
for all $j\ne i$). Then $\ka_0(\mu,\la)$ is the multiplicity of 
the part $\la_i$ in the partition $\la$.

In `exponential notation' for partitions,
$$
\mu=(1^{r_1(\mu)}2^{r_2(\mu)}\dots),\qquad
\la=(1^{r_1(\la)}2^{r_2(\la)}\dots),
$$
$\la\searrow\mu$ means that there exists a unique $k$ (equal to $\la_i$) 
such that 
$$
r_k(\la)=r_k(\mu)+1,\quad
r_{k-1}(\la)=r_{k-1}(\mu)-1,\quad
r_l(\la)=r_l(\mu) \quad\text{for $l\ne k$},
$$
and then
$$
\ka_0(\mu,\la)=r_k.
$$

The role of Young graph $\Y$ is played by the {\it Kingman graph\/}
$\K=(\Y,\ka_0)$: its vertices are the same as for $\Y$ (arbitrary
partitions) but each edge $\mu\nearrow\la$ turns into $\ka_o(\mu,\la)$
edges with the endpoints $\mu$ and $\nu$. The grading $\K=\cup\K_n$ of
the vertices remains unchanged. That is to say, the $n$th level $\K_n$
consists of partitions of $n$.

A new dimension function $\dim_0$ then arises: $\dim_0\la$ is still 
defined as the number of oriented paths from $\varnothing$ to $\la$ 
but we take into account edge multiplicities. The recurrence relation 
is modified as follows:
$$
\dim_0\la=\sum_{\mu:\,\mu\nearrow\la} \dim_0\mu\,\ka_0(\mu,\nu).
$$
By the very definition, the numbers $\ka_0(\mu,\la)$ are nothing but 
the coefficients in the expansion
$$
p_1^n=\sum_{\la:\in\K_n}\dim_0\la\cdot m_\la.
$$

The harmonicity condition is written as follows:
$$
\varphi(\mu)=\sum_{\la:\,\la\searrow\mu}\ka_0(\mu,\la)\,\varphi(\la).
$$

A {\it coherent system of distributions on the graph\/} $\K$ is a
sequence $M=(M_n)$ of probability distributions on the finite sets
$\K_n$ (partitions of $n$), subject to the coherence condition
$$
M_n(\mu)=\sum_{\la:\,\la\searrow\mu}
\frac{\dim_0\mu\cdot \ka_0(\mu,\la)}{\dim_0\la} \,M_{n+1}(\la),
\qquad \mu\in\K_n\,,
$$
equivalent to harmonicity of the function $M(\la)/\dim_0\la$ 
(here it is convenient to write $M(\mu)$ in place of $M_n(\mu)$). 
This concept is equivalent to that of {\it partition structure\/} 
in the sense of Kingman \cite{Ki1}.  

The next relation establishes a bijective correspondence 
$M\,\leftrightarrow\, P$ between coherent systems $M$ on Kingman's 
graph and probability measures $P$ on the simplex $\De$:
$$
\frac{M(\la)}{\dim_0\la}=\int_\De\wtm_\la(\al)\,P(d\al)\tag7.1
$$
(recall that $\wtm_\la$ is the image of $m_\la\in\La$ under the 
specialization $f\mapsto \wtf$ defined above). This result is 
{\it Kingman's theorem\/} \cite{Ki1, Ke1, KOO}. 

The coherent {\it $t$-systems\/} $M^{(t)}=(M^{(t)}_n)$ of 
distributions on Kingman's graph are defined as follows: $t>0$ 
is a parameter and
$$
M^{(t)}_n(\la)=\frac{t^{\ell(\la)}n!}{(t)_n\,  z_\la}, 
\qquad \la\in\K_n,\tag7.2
$$
where $\ell(\la)$ is the length   of the partition $\la$ (number 
of nonzero parts) and, in `exponential notation',
$$
z_\la=\prod_i\la_i\cdot\prod_k r_k(\la)!
=\prod_k k^{r_k(\la)}r_k(\la)!\,.\tag7.3
$$
The fact that $M^{(t)}$ is indeed a coherent system on $\K$  
is verified directly using the above formulas for $\dim_0\la$ 
and $\ka_0(\mu,\la)$. 

The systems $M^{(t)}$ are known as {\it Ewens partition structures\/} 
\cite{Ki1, Ki2}. According to Kingman's theorem, they determine certain 
probability measures on $\De$. The latter are called the 
{\it Poisson--Dirichlet distributions\/}, see \cite{Ki2}, and 
denoted as $PD(t)$. So, the link between $M^{(t)}$ and $PD(t)$ 
is as follows:
$$
M^{(t)}(\la)=\int_\De \wtm_\la(\al)(PD(t))(d\al).
$$

Given a probability measure $P$ on $\De$, we may regard it as a
measure on $\Om$. Hence, the definitions of the {\it control
measures\/} $\si_n$, of the {\it point process\/} $\PP$ attached to
$P$, and of the {\it correlation measures\/} $\rho_n$ of $\PP$ make
sense. As in the present situation there is no $\be$'s, the measure
$\si_n$ is concentrated on the cube $[0,1]^n$ contained in $[-1,1]^n$,
the state space of the process $\PP$ is the semiopen interval $(0,1]$,
and the $n$th correlation function $\rho_n$ lives on $(0,1]^n$.

Assume that $P$ is related to $M$ by \tht{7.1}. According to
\tht{3.1a}, the moments of $\si_n$ are given by
$$
\int_{[0,1]^n} x_1^{l_1}\dots x_n^{l_n}\,\si_n(dx_1\dots dx_n)=
\int_\De \wtp_{l_1+1}(\al)\dots \wtp_{l_n+1}(\al)\, P(d\al),
$$
which implies 
$$
\int_{[0,1]^n} x_1^{l_1}\dots x_n^{l_n}\,\si_n(dx_1\dots dx_n)
=\sum_\la [p_{l_1+1}\dots p_{l_n+1}:m_\la]
\,\frac{M(\la)}{\dim_0\la}\,. \tag7.4
$$
Here $\la$ ranges over $\K_{l_1+\dots l_n+n}$ and, for $f\in\La$, 
the symbol $[f\,:\,m_\la]$ denotes the coefficient of $m_\la$ in 
the expansion of $f$ into linear combination of monomial symmetric 
functions. It is worth noting that only partitions $\la$ with 
$\ell(\la)\le n$ really contribute to this formula, cf. \tht{3.1b}.

Consequently, given $M$, we have again a collection of moment problems
to determine the control measures $\si_n$, and from the measures
$\si_n$ one can get the correlation functions $\rho_n$. The procedure
is exactly the same as for the Young graph. We shall apply it 
to $M=M^{(t)}$. The corresponding control measures will be
denoted as $\sit_n$ and the correlation measures as $\rt_n$. The point
process attached to $PD(t)$ will be called the {\it Poisson--Dirichlet
process\/} and denoted as $\PDt$.

We shall see that in the present situation, calculation of the
correlation functions turns out to be much easier than in the case of
the Young graph. A formal explanations is that the coefficients in the
expansion of $p_{l_1+1}\dots p_{l_n+1}$ on monomial functions $m_\la$
are given by much simpler expressions than when expanding on Schur
functions $s_\la$.

To state the result we need some notation related to set partitions.
Recall (see \S4) that partitions of a set are the 
same thing as equivalence relations. We shall also need {\it
ordered\/} partitions, i.e., partitions with a fixed enumeration of
the blocs. As in Proposition 4.4, we denote the set of partitions of
$\{1,\dots,n\}$ with $r$ (nonempty) blocs by $\Pi(n,r)$; the set of
ordered partitions with $r$ blocs will be denoted as
$\wtPi(n,r)$. There is a natural projection $\wtPi(n,r)\to\Pi(n,r)$ ---
forgetting enumeration. For each $\pi\in\Pi(n,r)$, there are $r!$
ordered partitions $\wtpi$ over $\pi$.

Given $\pi\in\Pi(n,r)$, we define a `diagonal section' $[0,1]^n_\pi$ 
of the cube $[0,1]^n$ as intersection with all hyperplanes of the form 
$x_i=x_j$ where $i\sim j\mod\pi$. To an ordered partition 
$\wtpi=\wtpi_1\sqcup\dots\sqcup\wtpi_r$ over $\pi$ we assign a 
bijective map
$$
\gathered
[0,1]^r\,\to\,[0,1]^n_\pi\subset[0,1]^n\,,\\
y=(y_1,\dots,y_r)\,\mapsto\,x=(x_1,\dots,x_n),\\
\text{where $x_j=y_i$ provided that $j\in\wtpi_i$.}
\endgathered \tag7.5
$$
Define a measure 
$\whsit_{r,\wtpi}$ on the cube $[0,1]^r$ as follows:
$$
\whsit_{r,\wtpi}(dy)=t^r\cdot\prod_{i=1}^r y_i^{|\wtpi_i|-1}
\cdot(1-\sum_{i=1}^r y_i)^{t-1}_+\,dy,  \tag7.6
$$
and let $\sit_{n,r,\wtpi}$ be its image under the above map $y\mapsto
x$. The latter measure lives on $[0,1]^n_\pi$ and does not depend on
the choice of $\wtpi$ over $\pi$, whence we may denote is as
$\sit_{n,r,\pi}$.

\example{Example} For $n=3$ the list of partitions and measures 
is as follows.
\medskip 

$\bullet$ $\Pi(3,1)$:

$$
\pi=\{1,2,3\},
\quad \sit_{3,1,\pi}=t^3(1-x_1-x_2-x_3)^{t-1}_+dx
$$

$\bullet$ $\Pi(3,2)$:

$$
\gather
\pi=\{1\}\sqcup\{2,3\}, 
\quad \sit_{3,2,\pi}=t^2x_2(1-x_1-x_2)^{t-1}_+\delta(x_2-x_3)dx\\
\pi=\{2\}\sqcup\{3,1\}, 
\quad \sit_{3,2,\pi}=t^2x_3(1-x_2-x_3)^{t-1}_+\delta(x_3-x_1)dx\\
\pi=\{3\}\sqcup\{1,2\}, 
\quad \sit_{3,2,\pi}=t^2x_2(1-x_2-x_3)^{t-1}_+\delta(x_1-x_2)dx
\endgather
$$

$\bullet$ $\Pi(3,3)$:

$$
\pi=\{1\}\sqcup\{2\}\sqcup\{3\},
\quad \sit_{3,3,\pi}=t^3x_1^2(1-x_1)^{t-1}_+\de(x_1-x_2)\de(x_1-x_3)dx.
$$
\endexample

\proclaim{Theorem 7.1} In the above notation,
$$
\sit_n=\sum_{r=1}^n\;\sum_{\pi\in\Pi(n,r)}\sit_{n,r,\pi}\,.
$$
\endproclaim

\demo{Proof} We apply the general formula \tht{7.4} to $M=M^{(t)}$. 
Let us abbreviate
$$
|l|=l_1+\dots l_n\,, \qquad r_k=r_k(\la).
$$
Using the expressions for $M^{(t)}$ and $z_\la$ given in \tht{7.2} 
and \tht{7.3}, we get
$$
%\gather
\int_{[0,1]^n} x_1^{l_1}\dots x_n^{l_n}\,\sit_n(x)dx
=\sum_{r=1}^n\,\sum
\Sb \la:\,|\la|=|l|+n\\
\ell(\la)=r \endSb
[p_{l_1+1}\dots p_{l_n}:m_\la]\,
\frac{t^r\prod_{i=1}^r(\la_i-1)!}
{(t)_{|l|+n}\prod_{k\ge1}r_k!}\,. \tag7.7a
%\endgather
$$

Remark that for any $f\in\La$ and any $\la$ with $\ell(\la)=r$,
$$
[f:m_\la]=[f:x_1^{\la_1}\dots x_r^{\la_r}],
$$
where the square brackets on the right denote the coefficient of 
$x_1^{\la_1}\dots x_r^{\la_r}$ in the expansion of $f$ as linear 
combination of monomials. Therefore, we may rewrite \tht{7.7a} as 
$$
\multline
\int_{[0,1]^n} x_1^{l_1}\dots x_n^{l_n}\,\sit_n(x)dx\\
=\sum_{r=1}^n\,\sum
\Sb \la:\,|\la|=|l|+n\\
\ell(\la)=r \endSb
[p_{l_1+1}\dots p_{l_n}:x_1^{\la_1}\dots x_r^{\la_r}]\,
\frac{t^r\prod_{i=1}^r(\la_i-1)!}
{(t)_{|l|+n}\prod_{k\ge1}r_k!}\,. 
\endmultline  \tag7.7b
$$ 

Next, write $\la$ as the $n$-tuple $\la_1\ge\dots\ge\la_r>0$ and
remark that each summand in the right--hand side of \tht{7.7b} makes
sense for any $n$-tuple of positive integers $\la_1,\dots,\la_r$ and
is symmetric with respect to their permutations. Moreover, remark that
there are exactly $\frac{r!}{r_1!r_2!\dots}$ distinct permutations of
the numbers $\la_1,\dots,\la_r$. It follows that we may drop the
restriction $\la_1\ge\dots\ge\la_r$ and at the same time replace
$\prod r_k!$ by $r!$. Thus, we can transform \tht{7.7b} to 
 
$$
\multline
\int_{[0,1]^n} x_1^{l_1}\dots x_n^{l_n}\,\sit_n(x)dx\\
=\sum_{r=1}^n\,\sum
\Sb \la_1>0,\dots,\la_r>0\\
\sum\la_i=|l|+n \endSb
[p_{l_1+1}\dots p_{l_n+1}:x_1^{\la_1}\dots x_r^{\la_r}]\,
\frac{t^r\prod_{i=1}^r(\la_i-1)!}
{(t)_{|l|+n}\prod_{k\ge1}r_k!}\,.
\endmultline \tag7.8
$$

Write 
$$
p_{l_1+1}\dots p_{l_n+1}=(x_1^{l_1+1}+x_2^{l_1+1}+\dots)
\dots(x_1^{l_n+1}+x_2^{l_n+1}+\dots)
$$
and remove the parentheses. Then we will get a sum of monomials, each
of which corresponds to a choice of a summand from the first, second,
\dots, $n$th parentheses. We are interested only in monomials of the
form $x_1^{\la_1}\dots x_r^{\la_r}$ with strictly positive
$\la_1,\dots,\la_r$. There is a bijective correspondence between such
monomials and ordered partitions $\wtpi\in\wtPi(n,r)$. Specifically,
an index $j$ belongs to the $i$th bloc of $\pi$ if in the $j$th
parentheses, the $i$th summand was chosen. This yields a
correspondence
$$
\wtpi\mapsto\la=\la(\wtpi),\qquad
\la(\wtpi)_i=\sum_{j\in\wtpi_i}(l_j+1)
=(\sum_{j\in\wtpi_i}l_j)+|\wtpi_i|.
$$

Now we shall split the moment problem \tht{7.8} into a collection 
of moment problems corresponding to various $\wtpi$. Specifically, 
assume that for any $r=1,\dots,n$ and any $\wtpi\in\wtPi(n,r)$ we 
dispose of a measure $\sit_{n,r,\wtpi}$ on $[0,1]^n$ which solves 
the moment problem
$$
\int_{[0,1]^n} x_1^{l_1}\dots x_n^{l_n}\,\sit_{n,r,\wtpi}(x)\,dx
=\frac{t^r\prod_{i=1}^r(\la(\wtpi)_i-1)!}
{(t)_{|l|+n}}\tag7.9
$$
Then the measure
$$
\sit_n:=\sum_{r=1}^n\frac1{r!}
\sum_{\wtpi\in\wtPi(n,r)} \sit_{n,r,\wtpi} \tag7.10
$$
will solve the moment problem \tht{7.8}.

Looking at \tht{7.9} we remark that the $(l_1,\dots,l_n)$-moment
depends only of the sums
$$
m_i:=\sum_{j\in\wtpi_i}l_j\,,\qquad i=1,\dots,r.
$$
This indicates that the desired measure should live on the section
$[0,1]^n_\pi$, where $\pi$ stands for the (unordered) partition
corresponding to $\wtpi$. Then we identify $[0,1]^n_\pi$ with
$[0,1]^r$ via the map $y\mapsto x$ defined in \tht{7.5} and rewrite the
moment problem \tht{7.9} in terms of the coordinates $y_1,\dots,y_n$:
$$
\int_{[0,1]^r}y_1^{m_1}\dots y_r^{m_r}\,\widehat\si_r(y)\,dy=
\frac{t^r\prod_{i=1}^r(m_i+|\wtpi_i|-1)!}
{(t)_{\sum(m_i+|\wtpi_i|)}}\,, \tag7.11
$$
where $\widehat\si_r$ stands for the unknown measure on 
$[0,1]^r$ and we have used the identity
$$
|l|+n=\sum_{i=1}^r (m_i+|\wtpi_i|)
$$

On the other hand, consider the measure
$\widehat\si_{r,\wtpi}=\widehat\si^{(t)}_{r,\wtpi}$ as defined in
\tht{7.6}. This is a Dirichlet measure on a simplex (see
\cite{Ki2}) whose moments can be readily calculated (this is a
multivariate version of the classical Euler beta--integral). One
verifies that this measure solves the moment problem \tht{7.11}.

Thus, we have shown that the measure \tht{7.10} indeed solves the 
initial moment problem \tht{7.7}. Finally, we remark that 
the measure $\sit_{n,r,\wtpi}$ in \tht{7.10} actually depends 
on the image $\pi\in\Pi(n,r)$ of $\wtpi$, which concludes the 
proof. \qed
\enddemo

Theorem 7.1 yields new proofs of certain well--known properties of 
the Poisson--Dirichlet distributions $PD(t)$.

Consider the face $\De_0=\Om_0\cap\De$ of the simplex $\De$, i.e., 
$$
\De_0=\{(\al;\ga)\bigm | \ga=1-(\al_1+\al_2+\dots)=0\}.  
$$

\proclaim{Corollary 7.2 {(\rm cf. \cite{Ki, 9.4--9.5})}} 
The Poisson--Dirichlet distribution $PD(t)$ is concentrated on the
face $\De_0$ of $\De$.
\endproclaim

\demo{Proof} Applying Theorem 7.1 with $n=1$, we get 
$$
\sit_1(dx)=t(1-x)^{t-1}dx, \qquad 0\le x\le1.
$$
As this measure has no atom at zero, we conclude, by Proposition 3.4,
that $PD(t)$ lives on the face $\Om_0$ and hence on $\De_0$. \qed
\enddemo

\proclaim{Corollary 7.3} The Poisson--Dirichlet process $\PDt$ is
simple (see the definition at the beginning of \S4). 
\endproclaim

\demo{Proof} Apply Theorem 7.1 with $n=2$. Both $\Pi(2,2)$ and 
$\Pi(2,1)$ consist of a single element $\pi$: specifically,
$\pi=\{1\}\sqcup\{2\}$ and $\pi=\{1,2\}$, respectively. So, we have
$$
\sit_2=\sit_{2,2,\{1\}\sqcup\{2\}}+\sit_{2,1,\{1,2\}}.
$$ 
The first component on the right is an absolutely continuous
measure on the square $[0,1]^2$, while the second one is a singular
measure supported by the diagonal $x_1=x_2=y$ of $[0,1]^2$. According
to Proposition 3.5, we have to examine only the second component. It
is equal to 
$$
tx_1(1-x_1)_+^{t-1}\de(x_1-x_2)dx_1dx_2
$$ 
or, in terms of the coordinate $y$, to $ty(1-y)_+^{t-1}dy$. The latter
expression is just $\sit_1(dy)$ multiplied by $y$. {}From Proposition
3.5 we conclude that the process is simple. \qed
\enddemo

Note that this fact is evident from the construction of $PD(t)$ via a
subordinator, see \cite{Ki}.  

\proclaim{Corollary 7.4} The correlation
functions of the Poisson--Dirichlet process $\PDt$ are given by the formula 
$$
\rt_n(x_1,\dots,x_n)=
\frac{t^n(1-x_1-\dots-x_n)_+^{t-1}}
{x_1\dots x_n}\, \qquad n=1,2,\dots .  
$$ 
\endproclaim 

This result is due to Watterson \cite{W} but our approach differs from
that of \cite{W}. One more proof can be obtained by making use of the
fact that the so--called {\it size--biased sampling\/} from 
$PD(t)$ yields a sequence of {\it independent\/} random variables
\cite{Ki2, 9.6}.

\demo{Proof} Since the process is
simple (Corollary 7.3), all diagonal sections $[0,1]^n_\pi$ with
$\pi\ne\{1\}\sqcup\dots\sqcup\{n\}$ are negligible sets with respect
to the $n$th correlation measure $\rt_n$ (Proposition 4.3). Let
$$
([0,1]^n)'=[0,1]^n\cap(I^n)'=
\{x\in(0,1]^n\bigm| x_i\ne x_j, \quad i<j\}
$$
be the complement to all proper diagonal sections. By Theorem 7.1, all 
the components $\sit_{n,r,\pi}$ with $r<n$  are concentrated outside
$([0,1]^n)'$, and the only component with $r=n$ is absoltely
continuous with respect to $dx$ with density 
$t^n(1-x_1-\dots x_n)^{t-1}_+$. 

According to Proposition 4.3, the $n$th correlatiuon measure is obtained 
from the latter measure by dividing it by 
$|x_1\dots x_n|=x_1\dots x_n$. \qed
\enddemo

\example{Remark 7.6} Note that the structure of the $n$th controlling
measure, as described by Theorem 7.1, is in perfect accordance with
the decomposition given in Proposition 4.4.
\endexample

\head \S8. Appendix (A.~Borodin): a proof of theorem 2.1 \endhead

In this Appendix we present a simple direct proof of Theorem 2.1,
which is due to A.~Borodin.
 
We shall use the Frobenius notation
$\la=(p_1,\dots,p_d\,|\,q_1,\dots,q_d)$ for a Young diagram $\la$,
see \tht{2.5}. We start with the observation that the function
$\varphi:=\M/\dim\la$ can be written in the determinantal form
$$
\varphi(\la)=\frac{\det[m_{p_iq_j}]}{(t)_n}\,,
\qquad n=|\la|=\sum(p_i+q_i+1), \tag8.1
$$
where
$$
m_{p,q}=t\frac{(z+1)_p(z'+1)_p
(-z+1)_q(-z'+1)_q}{p!q!(p+q+1)};\quad t=zz'.\tag8.2
$$
Indeed, this easily follows from the formulas \tht{2.6}, \tht{2.7},
and the Cauchy formula
$$
\frac
{\prod\limits_{1\le i,j\le d}(p_i-p_j)(q_i-q_j)}
{\prod\limits_{1\le i\le d}\,
\prod\limits_{1\le j\le d}(p_i+q_j+1)}=
\det\left[\frac1{p_i+q_j+1}\right]\,.
$$

 Let us introduce a class of functions on the Young graph.
We fix a number $t\in \Bbb C\setminus \{0,-1,-2,\dots\}$, a sequence
$\{m_{pq}\}_{p,q=0}^{\infty}$ and set 
 $$
\varphi(\lambda)=\frac{\det[m_{p_i,q_j}]_{i,j=1}^d}{t(t+1)\cdots(t+n-1)}
\tag8.3
$$
\proclaim{Theorem 8.1} If a sequence $\{m_{pq}\}_{p,q=0}^{\infty}$ satisfies 
the relations$$
\gathered
 m_{p+1,q}+m_{p,q+1}-(p+q+1)m_{p,q}=m_{p,0}m_{0,q}\,,
\qquad p,q=0,1,\dots,\\
m_{0,0}=t
\endgathered\tag8.4
$$
then the function $\varphi$ defined by \tht{8.1} is harmonic. In other words,
$$
\varphi(\la)=\sum_{\nu\searrow\lambda}\varphi(\nu), \quad \la\in\Y.
\tag8.5
$$
\endproclaim

\demo{Proof} For a $l\times l$ matrix $A=(a_{ij})$ we shall denote by
$A\binom{i_1\dots i_k}{j_1\dots j_k}$ the determinant of the
submatrix of $A$ formed by the intersections of rows with numbers
$i_1\dots i_k$ and columns with numbers $j_1\dots j_k$. We shall also
denote by $A_{ij}$ the cofactor of $a_{ij}$. That is, 
$$
A_{ij}=(-1)^{i+j}A\binom{1\dots \hat i\dots l}{1\dots \hat j \dots l}.
$$
The transposed matrix to $(A_{ij})$ is equal to the inverse matrix of
$A$, multiplied by $\det A$, so that
$$
\sum_j a_{ij}A_{ij}=\sum_i a_{ij}A_{ij}=\det A.
$$
It follows that for any two sequences of numbers $v_1,\ldots,v_l$ and
$w_1,\ldots,w_l$ 
$$
\sum_{i,j=1}^l (v_i+w_j)a_{ij}A_{ij}
=\sum_{i=1}^l(v_i+w_i)\cdot\det A.\tag8.6 
$$ 

We proceed to verify the harmonicity relation \tht{8.5}. The Frobenius
coordinates of a diagram $\nu\searrow\la$ in \tht{8.5} are obtained from
the Frobenius coordinates of the diagram $\la$ by applying one of the
following three operations: 

1) $p_i\to p_i+1$ for a certain $i=1,\dots,d$, which corresponds to
creating a new box in the $i$th row above the diagonal;

2) $q_j\to q_j+1$ for a certain $j=1,\dots,d$, which corresponds to
creating a new box in the $j$th column below the diagonal;

3) adding a couple of coordinates $p_{d+1}=0$, $q_{d+1}=0$, which
corresponds to creating a new box on the diagonal.

It may happen that creating a new box in a certain position is
forbidden, because the resulting shape $\nu$ is not a Young diagram:
this occurs exactly when the set of the coordinates for $\nu$ 
contains two equal $p$-coordinates or two equal $q$-coordinates.
However, in such a case the formal application of formula \tht{8.3}
will give $\varphi(\nu)=0$ as the determinant in the numerator of \tht{8.3}
will vanish. This makes it possible to sum up over all the operations of
type 1), 2), 3), irrespective to whether the corresponding shape
$\nu$ is a Young diagram.  Then the harmonicity relation \tht{8.5} can be
rewritten in the following form (below we set $M=(m_{pq})$)
$$
\multline
\left(t+\sum_{i=1}^d(p_i+q_i+1)\right)M\binom{p_1\dots p_d}{q_1\dots q_d}
\\=\sum_{i=1}^d M\binom{p_1\dots p_i+1 \dots p_d}{q_1\dots q_d}
+\sum_{j=1}^d M\binom{p_1\dots p_d}{q_1\dots q_j+1\dots q_d}
+M\binom{p_1\dots p_d \ 0}{q_1\dots q_d \ 0}
\endmultline
$$

Expansion along the $i$th row gives
$$
M\binom{p_1\dots p_i+1 \dots p_d}{q_1\dots q_d}=\sum_{j=1}^d
(-1)^{i+j}m_{p_i+1,q_j}M\binom{p_1\dots \hat p_i \dots p_d}
{q_1\dots \hat q_j \dots q_d}.
$$
Similarly, expanding along the $j$th column, we get
$$
M\binom{p_1\dots  p_d}{q_1\dots q_j+1 \dots q_d}=\sum_{i=1}^d
(-1)^{i+j}m_{p_{i},q_j+1}M\binom{p_1\dots \hat p_i \dots p_d}
{q_1\dots \hat q_j \dots q_d}.
$$

Finally, expanding $M\binom{p_1\dots p_d \ 0}{q_1\dots q_d \ 0}$ along 
the last row and column we get
$$
\multline
M\binom{p_1\dots p_d \ 0}{q_1\dots q_d \ 0}\\ =
m_{0,0}M\binom{p_1\dots p_d}{q_1\dots q_d}+\sum_{i,j=1}^d
(-1)^{i+j+1}m_{p_i,0}m_{0,q_{j}}M\binom{p_1\dots \hat p_i \dots p_d}
{q_1\dots \hat q_j \dots q_d}
\endmultline
$$

Adding everything up, employing the assumption \tht{8.4} and applying
the relation \tht{8.6}, we get
$$
\align
&\sum_{i=1}^d M\binom{p_1\dots p_i+1 \dots p_d}{q_1\dots q_d}
+\sum_{j=1}^d M\binom{p_1\dots p_d}{q_1\dots q_j+1\dots q_d}
+M\binom{p_1\dots p_d \ 0}{q_1\dots q_d \ 0}\\
&=tM\binom{p_1\dots p_d}{q_1\dots q_d}\\&+\sum_{i,j=1}^d
(-1)^{i+j}(m_{p_i+1,q_j}+m_{p_i,q_j+1}-m_{p_i,0}m_{0,q_{j}})
M\binom{p_1\dots \hat p_i \dots p_d}
{q_1\dots \hat q_j \dots q_d}\\&=
tM\binom{p_1\dots p_d}{q_1\dots q_d}+\sum_{i,j=1}^d
(-1)^{i+j}(p_i+q_j)m_{p_i,q_j}M\binom{p_1\dots \hat p_i \dots p_d}
{q_1\dots \hat q_j \dots q_d}\\&=
\left(t+\sum_{i=1}^d(p_i+q_i+1)\right)
M\binom{p_1\dots p_d}{q_1\dots q_d}\,,
\endalign
$$
which concludes the proof. \qed
\enddemo

\proclaim{Corollary 8.2} The claim of Theorem 2.1 holds.
\endproclaim

\demo{Proof} It suffices to check that the sequence \tht{8.2} satisfies the
assumption \tht{8.4} of Theorem 8.1. But this is easily verified. \qed
\enddemo

\enddocument
\bye